\pdfoutput=1
\RequirePackage{ifpdf}
\ifpdf % We are running pdfTeX in pdf mode
\documentclass[pdftex]{sigma}
\else
\documentclass{sigma}
\fi

\usepackage[all]{xy}
\usepackage{enumerate}

\numberwithin{equation}{section}

\newtheorem{Theorem}{Theorem}[section]
\newtheorem{Corollary}[Theorem]{Corollary}
\newtheorem{Lemma}[Theorem]{Lemma}
\newtheorem{Proposition}[Theorem]{Proposition}
 { \theoremstyle{definition}
\newtheorem{Definition}[Theorem]{Definition}

\newtheorem{Example}[Theorem]{Example}
\newtheorem{Remark}[Theorem]{Remark} }

\begin{document}

\newcommand{\arXivNumber}{1110.0729}

\allowdisplaybreaks

\renewcommand{\PaperNumber}{052}

\FirstPageHeading

\ShortArticleName{Algebro-Geometric Solutions of the Generalized Virasoro Constraints}

\ArticleName{Algebro-Geometric Solutions \\ of the Generalized Virasoro Constraints}

\Author{Francisco Jos\'e PLAZA MART\'{I}N}

\AuthorNameForHeading{F.J.~Plaza Mart\'{\i}n}

\Address{Departamento de Matem\'aticas and IUFFYM, Universidad de Salamanca,\\  Plaza de la Merced 1-4,         37008 Salamanca, Spain}
\Email{\href{mailto:fplaza@usal.es}{fplaza@usal.es}}
\URLaddress{\url{http://mat.usal.es/~fplaza/}}

\ArticleDates{Received December 20, 2014, in f\/inal form July 02, 2015; Published online July 07, 2015}

\Abstract{We will describe algebro-geometric solutions of the KdV hierarchy whose $\tau$-functions in addition satisfy a generalization of the Virasoro constraints (and, in particular, a generalization of the string equation).
We show that these solutions are closely related to embeddings of the positive half of the Virasoro algebra into the Lie algebra of dif\/ferential operators on the circle. Our results are tested against the case of Witten--Kontsevich $\tau$-function. As by-products, we exhibit certain links of our methods with   double covers of the projective line equipped with a line bundle and with ${\rm Gl}(n)$-opers on the punctured disk.}


\Classification{14H71; 14D24; 81R10; 81T40}

\section{Introduction}

In recent decades, the relevant role in mathematics and in
physics of the class of functions that simultaneously satisfy the KdV-hierarchy
and the Virasoro constraints  has been discovered (as a~mere
illustration, let us
cite~\cite{DVV,DubrovinZhang,EynardOrantin,KacSchwarz,Kon}). On
the one hand, the KdV-hierarchy is currently quite well understood
after the works of Sato and it has been successfully applied in
topics such as Shiota's solution of the Schottky problem. On the other
hand, representation theory of Virasoro algebras has been
studied in depth; in particular its study has led to signif\/icant results in the theory of vertex operator algebras and in the geometric Langlands program. From the point of view of mathematical physics, these types of equations appear naturally in
conformal f\/ield theory when studying the symmetries of objects.

One principal character of this story is played by a particular
function in each theory. Among the variety of problems in which
both families of equations, KdV and Virasoro,  appear
together, let us pick up the functions of two of them: the
partition function of 2D gravity and the generating functions for
intersection numbers. The fact is that these functions are
$\tau$-functions for the KdV-hierarchy and solutions to the
Virasoro constraints. Indeed, this is an excellent example to
illustrate the deep consequences emerging from their study.
Recall, for instance, the Witten conjecture, proved by Kontsevich,
asserting that both $\tau$-functions do coincide or the
implications in the Gromov--Witten theory and enumerative geometry thanks to the ELSV formula,
which links Hurwitz numbers and Hodge integrals.

Within the context of matrix models,
Eynard and Orantin \cite{EynardOrantin,EynardOrantin2} have recently shown that
starting from an arbitrary spectral curve one can construct many
objects of the theory. In particular,  the nonperturbative
partition function is closely related to the theta function of the
curve and that, again, it is a $\tau$-function
for the KdV-hierarchy that solves the Virasoro constraints~\cite{EynardMarino}.

At this point we wondered about developing  a general formalism for the study of the KdV-Virasoro combination with emphasis on the role of the $\tau$-function. Here and henceforth, under \emph{Virasoro constraints for $\tau(t)$} we mean there exist a family of dif\/ferential operators $\{\bar L_k\}_{k\geq -1}$, satisfying the relations of the Virasoro algebra, such that $\bar L_k\tau(t)=0$ for $k\geq -1$. Thus, the Virasoro constraints in 2D quantum gravity as in~\cite{DVV,KacSchwarz,Kon} correspond to a particular choice of $\{\bar L_k\}_{k\geq -1}$. The present paper aims to of\/fer a f\/irst step in this program. The moral of the paper is that there are close and explicit relations among the following sets:
    \begin{gather}\label{eq:BigDiagramIntro}
    \begin{split}&
    \xymatrix@C=28pt@R=20pt{
    {\left\{\begin{gathered}
%        \text{functions in }
        \tau(t)\in {\mathbb C}[[t_1,t_2,\ldots]] \\
        \text{satisfying KdV and} \\
        \text{Virasoro constraints}
    \end{gathered}
    \right\} }
    \ar@<1ex>[r]^-A
    \ar@{-->}[d]^C
    &
    {\left\{\begin{gathered}
        \text{subalgebras } <\{\tilde L_n\}_{n\geq -1}> \\
        \text{ of }  {\mathcal D}^1_{{\mathbb C}((z))/{\mathbb C}}({\mathbb C}((z)),{\mathbb C}((z)))
        \\
        \text{with }[\tilde L_i,\tilde L_j]=(i-j)\tilde L_{i+j}
    \end{gathered} \right\} }
     \ar@{-->}[d]^D
    \ar@<1ex>@{-->}[l]^-B
    \\
    {\left\{\begin{gathered}
        \text{2:1-covers $X\to{\mathbb P}_1$ and a} \\
        \text{rank $1$ torsion free sheaf on }X
    \end{gathered} \right\} }
    &
    {\left\{\begin{gathered}
        \operatorname{Gl}(2)\text{-opers} \\
        \text{on }\operatorname{Spec}{\mathbb C}((z^{2}))
    \end{gathered} \right\} }
	}\end{split}
    \end{gather}
where we shall skip the details until later.

Let us discuss naively how the above arrows are constructed and postpone the rigor for the heart of the paper. Roughly, the map $A$ from functions to subalgebras is a forgetful map. More precisely, observe that saying that $\tau(t)$ is a $\tau$-function for the KdV that solves the Virasoro constraints means that $\tau(t)\in {\mathbb C}[[t_1,t_2,\ldots]]$ satisf\/ies the following equations:
    \begin{enumerate}[(i)]\itemsep=0pt
    \item KP-hierarchy,
    \item $\partial_{t_{2i}}\tau(t)=0$ (KdV hierarchy, provided that KP is fulf\/illed),
    \item $\bar L_k\tau(t)=0$, for $k\geq -1$ (Virasoro constraints), for certain dif\/ferential operators $\{\bar L_k\}_{k\geq -1}$ with $[\bar L_i,\bar L_j]=(i-j)\bar L_{i+j}$.
    \end{enumerate}
Recall that, in the case of 2D quantum gravity, the equation $\bar L_{-1}\tau(t)=0$ is the so called
\emph{string equation}. The following convention will used; $L_i$ will denote  a set of generators of the Witt algebra, as an abstract Lie algebra, such that $[L_i,L_j]=(i-j)L_{i+j}$; $\tilde L_i$ will be regarded as f\/ields acting on ${\mathbb C}((z))$ with the same Lie bracket; and, f\/inally, $\bar L_i$ will be f\/ields acting on the Fock space ${\mathbb C}[[t_1,t_2,\ldots]]$ with the same Lie bracket.

Now, the bosonization isomorphism  provides us with a method to
transport operators, $\bar L_k$, acting on ${\mathbb C}[[t_1,t_2,\ldots]]$
to  f\/irst-order dif\/ferential operators on ${\mathbb C}((z))$, which we
denote by $\tilde L_k$. Thus, $\{\tilde L_k\}_{k\geq -1}$ is the
desired subalgebra of the map $A$.

\looseness=-1
The dif\/f\/icult part is $B$, from subalgebras to functions,  since we need to construct a function. For this goal the Sato theory of soliton equations will be relied on heavily. Recall that the Sato theory
claims that $\tau$-functions for the KP hierarchy correspond
bijectively to certain subspaces of ${\mathbb C}((z))$ and that, moreover,
these subspaces are precisely the points of an inf\/inite-dimensional Grassmann manifold, $\operatorname{Gr}({\mathbb C}((z)))$. Hence, instead of looking for
$\tau$-functions for the KP we may look for subspaces $U\subseteq
{\mathbb C}((z))$ lying in the Grassmannian of ${\mathbb C}((z))$ and, at the end of
the day, the desired $\tau$-function will appear as the
$\tau$-function associated with the point $U$, $\tau_U(t)$. Thus, our strategy consists of restating the problem in terms of subspaces of ${\mathbb C}((z))$ with the help of Sato theory. Provided that the KP hierarchy is fulf\/illed, the  KdV hierarchy for $\tau(t)$ can be explicitly written in terms of subspaces; indeed, for a subspace $U$ it  is equivalent to the condition $z^{-2}U\subseteq U$.

So far, the problem of constructing the map $B$ can be restated as follows: given a Lie subalgebra $\langle\{\tilde L_k\}_{k\geq -1}\rangle$ of ${\mathcal D}^1_{{\mathbb C}((z))/{\mathbb C}}({\mathbb C}((z)),{\mathbb C}((z)))$ with $[\tilde L_i,\tilde L_j]=(i-j)\tilde L_{i+j}$, f\/ind $U\subset {\mathbb C}((z))$ such that the following three conditions are
satisf\/ied:
    \begin{enumerate}[(i$'$)]\itemsep=0pt
    \item $U\in\operatorname{Gr}({\mathbb C}((z)))$ (\emph{KP hierarchy}),
    \item $z^{-2}U\subset U$,
    \item $\tilde L_k U\subseteq U$ for $k\geq-1$,
    \end{enumerate}
which are the \emph{translations} of  items (i)--(iii) above.
Once this problem has been successfully solved, we associate with the subalgebra
the $\tau$-function of the subspace $U$, $\tau_U(t)$,  which, by construction, will  solve KdV and Virasoro simultaneously.

\looseness=1
Let us sketch the main steps of this construction. Our f\/irst task is
concerned with the study of Witt algebras. More explicitly, let
${\mathcal W}$ be the Witt algebra; that is, the ${\mathbb C}$-vector space
with basis $\{L_k\}_{k\in{\mathbb Z}}$ endowed with the Lie bracket
$[L_i,L_j]=(i-j) L_{i+j}$, and let ${\mathcal W}^+$ be the subalgebra
generated by $\{L_k\}_{k\geq -1}$. Then, Section~\ref{sec:representationWitt} is fully devoted to
the study of the maps of Lie algebras
\begin{gather*}
    \rho\colon \ {\mathcal W}^+\longrightarrow  {\mathcal D}^1,
\end{gather*}
where ${\mathcal D}^1:= {\mathcal D}^1_{{\mathbb C}((z))/{\mathbb C}}({\mathbb C}((z)),{\mathbb C}((z))))$ denotes the
algebra of f\/irst-order dif\/ferential operators of~${\mathbb C}((z))$. Then,
the Virasoro constraints will be written down in terms of the
operators $\tilde L_k:=\rho(L_k)$. Some results of this section
deserve special attention. For instance,
Theorem~\ref{thm:VirasoroSubAlg} deals with explicit expressions
for~$\rho$; it claims that there is a 1-1-correspondence:
    \begin{gather}\label{eq:1-1intro}
    \left\{\begin{matrix}
    \rho\in\operatorname{Hom}_{\text{Lie-alg}}\big({\mathcal
    W}^+,{\mathcal D}^1\big)
   \\
    \text{ s.t. }\rho\neq 0
    \end{matrix}\right\}
    \, \overset{\text{1-1}}\longleftrightarrow \, % \ar@{<-->}[r]^{1-1} &
    \left\{\begin{matrix}
    (h(z),c,b(z))\text{ s.t. }h'(z)\in{\mathbb C}((z))^*
    \\
    c\in{\mathbb C},\; b(z)\in{\mathbb C}((z))
    \end{matrix}\right\},
    \end{gather}
which is explicitly given by
\begin{gather*}
    \rho(L_i) =
    \frac{-h(z)^{i+1}}{h'(z)} \partial_z  -
    (i+1) c \cdot h(z)^i + \frac{h(z)^{i+1}}{h'(z)}b(z).
\end{gather*}
Other results are concerned with the classif\/ication of
$\operatorname{Im}\rho$ as a subalgebra of ${\mathcal D}^1$ (Theorem~\ref{thm:ImRho1=ImRho2}), with its  conjugacy class (Theorem~\ref{thm:HomModAutxC((z))})
as well as with the isomorphism classes of its central
extensions (Theorem~\ref{thm:IsomClassVirasoro}).

The second main problem to be tackled is the explicit description
of the vector space $U$, which is the main topic in
Section~\ref{sec:SolutionsTAUVir&KP}. To begin with, we examine what the
stabilizer of a subspace stable under $\rho$ looks like, and,
since it turns out to be ${\mathbb C}[h(z)]$
(Theorem~\ref{thm:stab=C[h]orC}), we conclude that, for the goals
of this paper, we may restrict our study to the case of those maps
$\rho$ with $h(z)=z^{-2}$ under the correspondence~(\ref{eq:1-1intro}). Now, let us comment on a technical issue
which is crucial in our approach. Observe that for any
$1$-dimensional ${\mathbb C}((z))$-vector space $V$, it holds that the
bosonization isomorphism establishes an identif\/ication
$H^0(\operatorname{Gr}(V),\operatorname{Det}^*)=\Lambda^{\frac{\infty}{2}}V
\simeq {\mathbb C}[[t_1,t_2,\ldots]]$ that is compatible with the actions
of the Virasoro algebra on $\Lambda^{\frac{\infty}{2}}V$ and that
of the Witt algebra on $\operatorname{Gr}(V)$. Thus, instead of ${\mathbb C}((z))$ we may
use any $1$-dimensional ${\mathbb C}((z))$-vector space~$V$. Note that
replacing one vector space by another amounts to conjugating the
operators. Regarding $\tau$-functions, which are our object of
interest, observe that two subpaces $U_1\in \operatorname{Gr}(V_1)$ and $U_2\in
\operatorname{Gr}(V_2)$ will share the same $\tau$-function as long as there is
an isomorphism $V_1\overset{\sim}\to V_2$ mapping~$U_1$ to~$U_2$.
An explanation for this lies in the property that the
$\tau$-function of a point $U$  \emph{depends only on the
coordinates} of vectors of $U$ w.r.t.\ a~basis of $V$ and not on
$V$ or $U$ themselves. Thus, $({\mathbb C}((z)),\rho)$ can be replaced by another $1$-dimensional ${\mathbb C}((z))$-vector space endowed with the action induced by $\rho$ by conjugation.

Then, it will be proved  that, for each solution $w(z)$ of
the Airy equation, there is a pair $(V^{wv},\rho^{wv})$ as above, where
$v$ depends only on $\rho$,  and a subspace ${\mathcal U}(w)$ of
$V^{wv}$ satisfying ${\mathcal U}(w)\in \operatorname{Gr}(V^{wv})$, $z^{-2}{\mathcal
U}(w)\subset {\mathcal U}(w)$ and $\rho^{wv}(L_k){\mathcal
U}(w)\subseteq {\mathcal U}(w)$. It follows that $\tau_{ {\mathcal
U}(w)}(t)\in{\mathbb C}[[t_1,t_3,\ldots]]$ is the desired function that
solves KdV and Virasoro simultaneously. This is the main result of
the paper (Theorem~\ref{thm:ExistenceU(v)}).

Once the function has been found, we may wonder about the
dependence on the choice of~$w$ and whether there is more than one
map producing the same function. The answer is that
$\partial_{t_{2i+1}}\partial_{t_{2j+1}}\log \tau_{
{\mathcal U}(w)}(t)$, for $i,j\geq 0$, does not depend on the
$w$ chosen (Theorem~\ref{thm:indep}). Further, the following result on uniqueness holds: let $\rho_1$ and $\rho_2$ be given and let $\tau_1(t)$ and $\tau_2(t)$ be two solutions of the KdV hierarchy and of the Virasoro constraints attached to~$\rho_1$ and to~$\rho_2$; then, the second derivatives of their logarithms coincide if and only if the Virasoro constraints for~$\rho_1$ and~$\rho_2$ coincide (Theorem~\ref{thm:uniqueness}). In a simpler way, the \emph{theory} (i.e., KdV plus Virasoro) determines uniquely the $\tau$-function and conversely.

\looseness=1
In the last section (see~Section~\ref{sec:ExamplesVirasoro}) some applications are of\/fered. First, our results are \emph{tested} against the well known case of 2D gravity.  Using~\cite{DVV,Douglas,Givental,KacSchwarz,Kaz,Kon} as the guidelines
of the theory, we shall check that our formalism f\/its perfectly within the algebraic manipulation of the Virasoro algebra appearing in the 2D gravity. In particular, the fact that the
partition function satisf\/ies KdV and Virasoro follows from the
very construction and does not need further proof.

A salient by-product emerge from a closer analysis of the
subspace $U$ and whose study deserves future research.  First in Section~\ref{subsec:spectral},  Krichever theory allows us to give an alternative description of~$U$ in terms of algebraic curves. Namely, under the hypotheses of Theorem~\ref{thm:ExistenceU(v)}, it holds that~$U$ is def\/ined by a 2:1-cover $X\to{\mathbb P}_1$ together with a rank $1$ torsion free sheaf on~$X$.
Moreover, certain variations of the curve would be governed by the Painlev\'{e} equations and correspond to isomonodromic deformations. Second, one can easily check that the pair $({\mathbb C}((z)),\rho)$ def\/ines a $\operatorname{Gl}(n)$-oper on the punctured disk, $\operatorname{Spec}{\mathbb C}((h(z)^{-1}))$. We hope that our construction might be better understood within the framework of Frenkel's approach to the geometric Langlands program~\cite{FrLoop}.

{\sloppy On the other hand,  the pair $({\mathbb C}((z)),\rho)$ def\/ines a $\operatorname{Gl}(n)$-oper on the punctured disk, $\operatorname{Spec}{\mathbb C}((h(z)^{-1}))$. Section~\ref{subsec:oper} is concerned with this construction which is the $D$ arrow of the above diagram for $h(z)=z^{-2}$. It would be natural to  expect a generalization of the relation between actions of a formal loop group (see Corollary~\ref{cor:HomLiesl2D1}) and opers on the punctured disk and, therefore,  we believe that our construction might be better understood within the framework of Frenkel's approach to the geometric Langlands program~\cite{FrLoop}.

}

Finally, we show that the families of Virasoro algebras used in~\cite{Kaz,LiuXu,MirzakhaniInvent,MirzakhaniJAMS,MulaseSafnuk,MulaseZhang} for the study of the topological recursion also f\/it into our framework and have a natural geometrical interpretation (see Section~\ref{subsec:UniverFamily}). In particular, it's shown that a family of $\tau$-functions satisfying KdV and Virasoro is equivalent to a family of actions of the Witt algebra. It is worth noticing that there are instances of such $1$-parameter families (e.g., the sine curve in~\cite{MulaseSafnuk}) which are indeed the spectral curve in Eynard--Orantin theory.

Before f\/inishing this introduction, let us mention some problems
that, in our opinion, can be studied with our techniques. For
instance, regarding  Theorem~\ref{thm:uniqueness},
we
point out that these functions arise in the study of monodromy-preserving deformations as well as in the theory of Frobenius
manifolds (e.g., \cite{DubrovinZhang} and \cite[Section~8]{Givental}). Second, the families of actions given
in~Section~\ref{subsec:UniverFamily}  suggests the possibility of
connecting our approach with the computation of the $\tau$-function
as a matrix integral (e.g.,~\cite{AMSvM}).
Finally, we shall extend our formalism to include special second-order dif\/ferential operators in order to study \emph{loop equations} and \emph{cut-and-join operators} and to consider the case of the $n$-KdV hierarchy.

\section{Actions of Witt algebras}\label{sec:representationWitt}

Let ${\mathcal W}$ be the Witt algebra; that is, the Lie algebra that is freely generated by $\{L_k\,\vert\, k\in {\mathbb Z}\}$ as a~${\mathbb C}$-vector space and endowed with the following Lie bracket  $[ L_i, L_j ]  =  (i-j) L_{i+j}$. Let us denote by ${\mathcal W}^+\subset {\mathcal W}$ the Lie subalgebra generated by $\{L_k\vert k\geq -1\}$.

Let $V$ denote a $1$-dimensional ${\mathbb C}((z))$-vector space and let ${\mathcal D}^1_{{\mathbb C}((z))/{\mathbb C}}(V,V)$ be the Lie algebra generated by f\/irst-order dif\/ferential operators. The symbol map is{\samepage
\begin{gather*}
    \sigma \colon \  {\mathcal D}^1_{{\mathbb C}((z))/{\mathbb C}}(V,V)\longrightarrow \operatorname{Der}_{{\mathbb C}}\big({\mathbb C}((z))\big)
    = {\mathbb C}((z))\partial_{z} \overset{\sim}\to {\mathbb C}((z)),
\end{gather*}
where the last map sends $f(z)\partial_z$ to $f(z)$.}

\subsection{General form}\label{subsec:ComputationVirasoroSubalg}

We are interested in pairs $(V,\rho)$ consisting of a $1$-dimensional ${\mathbb C}((z))$-vector space, $V$, and a Lie algebra homomorphism:
\begin{gather*}
    {\mathcal W}^+
     \overset{\rho}\longrightarrow
    {\mathcal D}^1_{{\mathbb C}((z))/{\mathbb C}}(V,V).
\end{gather*}
For the sake of brevity and when no confusion arises, the pair
$(V,\rho)$ will be called an {\it action} of~${\mathcal W}^+$.

In this subsection we shall focus on the case of $V={\mathbb C}((z))$.
Nevertheless, similar results can be proved for $V$ arbitrary by
f\/ixing an isomorphism $V\simeq{\mathbb C}((z))$ (see
Section~\ref{subsec:conjugation} for the dependence on the choice of the
isomorphism).

Although a lot of properties on ${\mathcal W}$-modules are known, the following explicit description of its actions on ${\mathcal D}^1_{{\mathbb C}((z))/{\mathbb C}}({\mathbb C}((z)),{\mathbb C}((z)))$ seems to be brand new (see also Remark~\ref{rem:densities}) except for the result of Kirillov when studying discrete series of unirreps for Virasoro in terms of univalent functions~\cite{Kirillov}.

\begin{Theorem}\label{thm:VirasoroSubAlg}
Let $V$ be ${\mathbb C}((z))$ and $\rho\colon {\mathcal W}^+\to
{\mathcal D}^1_{{\mathbb C}((z))/{\mathbb C}}(V,V)$ be a ${\mathbb C}$-linear map such that $\rho\neq 0$.
Then, the map $\rho$ is a Lie
algebra homomorphism if and only if there exist functions $h(z)
, b(z) \in {\mathbb C}((z))$ and a constant $c\in {\mathbb C}$ such that
$h'(z)\neq 0$ and
    \begin{gather}\label{eq:rhoLk-hbc}
    \rho(L_i)  =
    \frac{-h(z)^{i+1}}{h'(z)} \partial_z  -
    (i+1) c   h(z)^i + \frac{h(z)^{i+1}}{h'(z)}b(z).
    \end{gather}
\end{Theorem}

\begin{proof}
The converse is a straightforward computation. Let us prove the direct one.

Let us write $\rho(L_k) = a_k(z) \partial_z  + b_k(z)$. Since $\rho$ is a map of Lie algebras, the expression of the bracket $[ L_i, L_j ]   =   (i-j) L_{i+j}$ implies the following set of equations
    \begin{gather}
      a_i(z) a_j'(z) - a_j(z) a_i'(z) = (i-j) a_{i+j}(z), \label{eq:diffeqforcoeffa}
    \\
      a_i(z) b_j'(z) - a_j(z) b_i'(z) = (i-j) b_{i+j}(z).\label{eq:diffeqforcoeffb}
    \end{gather}

Observe that if $a_{-1}(z)=\sigma(\rho(L_{-1}))= 0$, then we let $i$ be equal to $-1$ in equation (\ref{eq:diffeqforcoeffa}) and have that $a_j(z)=0$ for all $j\geq -1$. Substituting in equation (\ref{eq:diffeqforcoeffb}), it follows that $\rho\equiv 0$.

Hence, we now assume that $a_{-1}(z)=\sigma(\rho(L_{-1}))\neq 0$. Let us f\/ix $L_{-1}$ and solve this system in terms of its
coef\/f\/icients.

Letting $i=-1$ in equation~(\ref{eq:diffeqforcoeffa}), dividing  by
$a_{-1}(z)^2$ and integrating, it follows that $a_j(z)=
-(1+j)a_{-1}(z)\int^z a_{j-1}(t)a_{-1}(t)^{-2} \operatorname{d}t$. Hence,  $a_j(z)$
can be determined recursively from $a_{-1}(z)$. Indeed, the case
$j=0$ yields $a_0(z)=a_{-1}(z)(\alpha-\int^z
\frac{\operatorname{d}t}{a_{-1}(t)})$ for $\alpha\in{\mathbb C}$. Since
$a_{-1}(z),a_0(z)\in{\mathbb C}((z))$, it follows that
$(\alpha-\int^z \frac{\operatorname{d}t}{a_{-1}(t)})$ must lie in
${\mathbb C}((z))$; i.e., there exists $h(z)\in {\mathbb C}((z))$
such that
\begin{gather*}
    a_{-1}(z) = \frac{-1}{h'(z)}.
\end{gather*}
Thus, setting the free term of $h(z)$ to be equal to that constant, it
follows that
\begin{gather*}
    a_0(z) = \frac{-h(z)}{h'(z)}.
\end{gather*}
Now, induction procedure proves straightforwardly that
\begin{gather*}
    a_{i}(z) = \frac{-h(z)^{i+1}}{h'(z)}.
\end{gather*}
Let us now focus on the $b_i$ variables. Firstly, let us deal with
the case $h(z)=z^n$; hence,  $a_i(z)=-\frac{1}{n}z^{n i +1}$ and
equation~(\ref{eq:diffeqforcoeffb}) has the following shape
    \begin{gather}\label{eq:diffeqforcoeffbh=z^n}
    -\frac{1}{n}z^{n i +1} b_j'(z) + \frac{1}{n}z^{n j +1} b_i'(z) = (i-j) b_{i+j}(z).
    \end{gather}
Let us write $b_j(z)$ as $\sum_{k} b_{j,k} z^k$, where $b_{j,k}=0$
for $k\ll 0 $. Computing the coef\/f\/icients of $z^k$, in
equation~(\ref{eq:diffeqforcoeffb}) one has the relation
\begin{gather*}
    \left( -\frac1n (k-n i) b_{j,k-n i} + \frac1n(k-n j) b_{i,k-n j}\right)  =  (i-j) b_{i+j,k}.
\end{gather*}
The case $j=0$ implies that $(k-n i)b_{i,k} = (k- n i) b_{0,k- n
i}$ and, therefore $b_{i,k} =  b_{0,k- n i}$ for $k \neq n i$;
that is, the dif\/ference between $z^{-ni}b_i(z)$ and $b_0(z)$ is a
constant.  Expressing this condition in terms of $b_{-1}(z)$,  the
following formula for $b_i(z)$ holds
    \begin{gather}\label{eq:bi1}
    b_i(z)   =  \big(c_i + b_{-1}(z)z^n\big) z^{n i}
    \end{gather}
for some $c_i\in {\mathbb C}$ and $c_{-1}=0$. Plugging this into
equation~(\ref{eq:diffeqforcoeffb}) and setting $i$ equal to $-1$,
we f\/ind a constraint for  the $c_i$
\begin{gather*}
    j c_j + c_{-1}- (j+1)c_{j-1}  =  0, \qquad c_{-1}=0,
\end{gather*}
whose general solution is
    \begin{gather}\label{eq:lambdaj}
    c_{j}  =  -c\cdot  (j+1)
    \end{gather}
for a complex number $c=-c_0\in{\mathbb C}$. Bearing in mind that $h'(z)$ is
invertible, there is no harm in assuming that $b_{-1}(z)$ is of
the form $\frac1{h'(z)}b(z)$.  Thus, from equations~(\ref{eq:bi1})
and~(\ref{eq:lambdaj}), the general solution for the case
$h(z)=z^n$ is
    \begin{gather}\label{eq:SOLdiffeqforcoeffbh=z^n}
    b_i(z)   =  \left(-(i+1)c+ z^n \frac{b(z)}{nz^{n-1}}\right) z^{n i}.
    \end{gather}

The general case, i.e., for $h(z)$ arbitrary, follows from the fact
that there is a ${\mathbb C}$-algebra automorphism of ${\mathbb C}[[z]]$, $\phi$,
such that $\phi(h(z))=z^n$ where $h(z)=a_n
z^n+a_{n+1}z^{n+1}+\cdots$ and $a_n\neq 0$. That is, in order to
solve equation~(\ref{eq:diffeqforcoeffb}), we consider $\phi$,
such that $\phi(h(z))=z^n$. We transform
equation~(\ref{eq:diffeqforcoeffb}) by $\phi$, which is
equation~(\ref{eq:diffeqforcoeffbh=z^n}), and consider its
solutions~(\ref{eq:SOLdiffeqforcoeffbh=z^n}). Thus, transforming
the solutions by the inverse automorphism, $\phi^{-1}$, we have
that the general solution for equation~(\ref{eq:diffeqforcoeffb})
is as follows
\begin{gather*}
    b_i(z)
  =  \left(-(i+1)c+ h(z) \frac{b(z)}{h'(z)}\right) h(z)^i.\tag*{\qed}
\end{gather*}
\renewcommand{\qed}{}
\end{proof}

Observing that in the previous proof only the operators
$\rho(L_i)$ for $i=-1,0$ were necessary, we have the following:

\begin{Corollary}\label{cor:HomLiesl2D1}
The  restriction maps
\begin{gather*}
    {\mathfrak{sl}}_2({\mathbb C})\simeq \langle L_{-1},L_0,L_1\rangle   \subset   {\mathcal W}^+
  \subset    {\mathcal W}
\end{gather*}
induce isomorphisms
\begin{gather*}
    \operatorname{Hom}_{\text{\rm Lie-alg}}\big({\mathcal W}, {\mathcal D}^1\big)\, \overset{\sim}\to \,
    \operatorname{Hom}_{\text{\rm Lie-alg}}\big({\mathcal W}^+, {\mathcal D}^1\big)\, \overset{\sim}\to \,
    \operatorname{Hom}_{\text{\rm Lie-alg}}\big({\mathfrak{sl}}_2({\mathbb C}), {\mathcal D}^1\big),
\end{gather*}
where ${\mathcal D}^1$ denotes ${\mathcal D}^1_{{\mathbb C}((z))/{\mathbb C}}({\mathbb C}((z)),{\mathbb C}((z)))$. Under these isomorphisms, the ${\mathfrak{sl}}_2({\mathbb C})$ representation associated to  $\rho$ $($attached to $(h(z),c,b(z)))$ has Casimir operator $1+4c^2$.
\end{Corollary}

\begin{Remark}
The combination of Corollary~\ref{cor:HomLiesl2D1} and formula~(\ref{eq:rhoLk-hbc}) with $h(z)=z^{-1}+d$ for $d\in{\mathbb C}$ agrees with Mulase's result on injective Lie algebra morphisms from ${\mathfrak{sl}}_2({\mathbb C})$ to ${\mathcal D}^1$ \cite[Proposition~8.6]{Mulase-AlgKP}.
\end{Remark}

\begin{Remark}The reason for restricting
ourselves to the case of f\/irst-order dif\/ferential operators of $V$
instead of more general operators relies on the fact that this is
the case in most situations, e.g., 2D gravity. From a more
theoretical point of view, the equivalence of categories between
Atiyah algebras and dif\/ferential operator algebras
\cite{BeilinsonSchehtman}, means that ${\mathcal D}^1_{{\mathbb C}((z))/{\mathbb C}}(V,V)$
is the natural candidate to begin with. Consequently, one expect a link with he moduli space of pairs consisting of a curve and a line
bundle on it, since the Lie group associated with
${\mathcal D}^1_{{\mathbb C}((z))/{\mathbb C}}(V,V)$ uniformizes this moduli space~\cite{HGPM-sgl}. By studying higher order dif\/ferential operators, one may relate this with the study  $W$-algebras (and their
representation theory, etc.~\cite{FKRW}, \cite[Section~1.2]{AMSvM}).\end{Remark}

\begin{Remark}\label{rem:LieBraLkh^j}
It holds that
    \begin{gather}\label{eq:LieBraLkh^j}
    \big[ \rho(L_k) , h(z)^j \big]  =  -j \cdot h(z)^{k+j}
    \end{gather}
as ${\mathbb C}$-linear operators on $V$. For  $g(z)\in{\mathbb C}((z))$ arbitrary,
it holds that $[\rho( L_k) , g] = -\frac{h(z)^{k+1}}{h'(z)}
g'(z)$.
\end{Remark}

\begin{Corollary}\label{cor:RhoSigma}
If $\sigma(\rho(L_{-1}))\neq 0$, then $\rho$ and $\sigma\circ\rho$ are injective and $\operatorname{Im}(\sigma\circ\rho)=\frac{1}{h'(z)}{\mathbb C}[h(z)]$.
\end{Corollary}

\begin{proof}
Let us prove that $\sigma\circ\rho$ is injective. Let us assume
that a linear combination $\sum_k \lambda_kL_k$ lies on the kernel
of $\sigma\circ\rho$; that is,  $\sigma(\rho(\sum_k
\lambda_kL_k))=0$. The previous Theorem yields the identity
$\sigma(\rho( L_k )) = -\frac{h(z)^{k+1}}{h'(z)}$, and therefore
we have that
\begin{gather*}
    0  =  \sigma\left(\rho\left(\sum_k \lambda_kL_k\right)\right)
     =
    \left( \sum_k \lambda_k h(z)^{k+1} \right) \sigma\big(\rho(L_{-1})\big).
\end{gather*}
Since $\sigma(\rho(L_{-1}))=-\frac1{h'(z)}\neq 0$, it follows that $\lambda_k=0$ for all $k$, and the claim is proved.

Second, the injectivity of $\sigma\circ\rho$ implies the injectivity of $\rho$. Finally, $\operatorname{Im}(\sigma\circ\rho)=\frac{1}{h'(z)}{\mathbb C}[h(z)]$ follows easily from equation~(\ref{eq:rhoLk-hbc}).
\end{proof}

Let ${\mathfrak v}\colon{\mathbb C}((z))\to {\mathbb Z}\cup\{\infty\}$ be the valuation
associated with $z$; that is, ${\mathfrak v}(0)=\infty$ and, for
$h(z)\neq 0$,  ${\mathfrak v}(h(z))=a$ if\/f $a$ is the largest
integer number such that $h(z)\in z^a{\mathbb C}[[z]]$ and, in this
situation, $a$ will be called the {\it order} of $h$.

\begin{Theorem}\label{thm:ImRho1=ImRho2}
For $i=1,2$, let $h_i(z), b_i(z)\in {\mathbb C}((z))$ and $c_i\in{\mathbb C}$ with
$h'_i(z)\neq 0$. Let $(V={\mathbb C}((z)),\rho_i)$ be the action of
${\mathcal W}^+$ associated with elements $h_i(z)$, $b_i(z)$, $c_i$,
as in Theorem~{\rm \ref{thm:VirasoroSubAlg}}.

If  $\operatorname{Im} \rho_1=\operatorname{Im}\rho_2$ and the signs of ${\mathfrak v}(h_1(z))$ and ${\mathfrak v}(h_2(z))$ are equal and negative, then $b_1(z)=b_2(z)$, $c_1=c_2$ and  $h_1(z)=\alpha h_2(z)+\beta$ for some $\alpha\in{\mathbb C}^*$, $\beta\in{\mathbb C}$.

Conversely, if $b_1(z)=b_2(z)$, $c_1=c_2$ and  $h_1(z)=\alpha h_2(z)+\beta$ for some $\alpha\in{\mathbb C}^*$, $\beta\in{\mathbb C}$, then $\operatorname{Im} \rho_1=\operatorname{Im}\rho_2$. Moreover, there exists a Lie algebra automorphism $\phi$ of ${\mathcal W}^+$ such that $\rho_2 =\rho_1 \circ\phi$.
\end{Theorem}

\begin{proof}
From the hypothesis $\operatorname{Im} \rho_1=\operatorname{Im}\rho_2$, it holds that there exist $\{\lambda_{kl}\,\vert\,  k,l\geq -1\}$ such that
\begin{gather*}
    \rho_1(L_k) =  \sum_{l\geq -1} \lambda_{kl} \rho_2(L_l).
\end{gather*}
By the explicit expression obtained in Theorem~\ref{thm:VirasoroSubAlg}, this identity is equivalent to the equations
    \begin{gather}\label{eq:ThmEquiSymbol}
    \frac{h_1(z)^{k+1}}{h_1'(z)}  =
    \sum_{l\geq -1} \lambda_{kl} \frac{h_2(z)^{l+1}}{h_2'(z)}
    \end{gather}
and
    \begin{gather}\label{eq:ThmEquiEnd}
     \frac{h_1(z)^{k+1}}{h_1'(z)} b_1(z) - (k+1) c_1 h_1(z)^k
   =
    \sum_{l\geq -1} \lambda_{kl} \left( \frac{h_2(z)^{l+1}}{h_2'(z)} b_2(z) - (l+1) c_2 h_2(z)^l\right).
    \end{gather}
Observe that the derivative of equation~(\ref{eq:ThmEquiSymbol}) w.r.t.~$z$ yields
    \begin{gather}
    (k+1) h_1(z)^k - \frac{h_1(z)^{k+1} h_1''(z)}{h_1'(z)^2}
 =
    \sum_{l\geq -1} \lambda_{kl} \left((l+1) h_2(z)^l - \frac{h_2(z)^{l+1} h_2''(z)}{h_2'(z)^2}\right)
  \nonumber \\
\hphantom{(k+1) h_1(z)^k - \frac{h_1(z)^{k+1} h_1''(z)}{h_1'(z)^2}}{}
   =
    \sum_{l\geq -1} \lambda_{kl} (l+1) h_2(z)^l  -   \frac{h_1(z)^{k+1}}{h_1'(z)}\cdot \frac{ h_2''(z)}{h_2'(z)}.\label{eq:ThmEquiSymbolDerivative}
    \end{gather}
One computes equation~(\ref{eq:ThmEquiEnd}) plus
equation~(\ref{eq:ThmEquiSymbol}) times $(-b_2(z))$ plus
equation~(\ref{eq:ThmEquiSymbolDerivative}) multiplied by
$c_2$, and one obtains
\begin{gather*}
    \frac{h_1(z)^{k+1}}{h_1'(z)}\left( (b_1(z)-b_2(z)) - (k+1)(c_1-c_2) \frac{h_1'(z)}{h_1(z)} -
    \left(\frac{h_1''(z)}{h_1'(z)} - \frac{h_2''(z)}{h_2'(z)}\right) c_2 \right)  =  0.
\end{gather*}
Since this holds for all $k\geq -1$, it follows that
    \begin{gather}
    c_1-c_2 = 0,
\qquad
    (b_1(z)-b_2(z)) - \left(\frac{h_1''(z)}{h_1'(z)} - \frac{h_2''(z)}{h_2'(z)}\right) c_2  =  0.
    \label{eq:ThmEquiEnd2}
    \end{gather}
Hence $c_1=c_2$.

Further, Corollary~\ref{cor:RhoSigma} shows that
$\frac{1}{h_1'(z)}{\mathbb C}[h_1(z)] = \frac{1}{h'_2(z)}{\mathbb C}[h_2(z)]$.
Hence, there are polyno\-mials~$p_i$ such that $\frac{1}{h_1'(z)}=
\frac{p_2(h_2(z))}{h'_2(z)}$ and $\frac{1}{h_2'(z)} =
\frac{p_1(h_1(z))}{h'_1(z)}$. These identities imply that
\begin{gather*}
    p_1(h_1(z))  p_2(h_2(z))  =  1
\end{gather*}
and, thus
\begin{gather*}
    \operatorname{deg}(p_1){\mathfrak v}(h_1(z)) +
    \operatorname{deg}(p_2){\mathfrak v}(h_2(z))   =  0.
\end{gather*}
The assumption about the signs of ${\mathfrak v}(h_i(z))$ implies
that $p_i$ is constant for $i=1,2$, say $p_1(x)=\alpha\in{\mathbb C}^*$.
And, therefore, $h_1'(z)=\alpha h_2'(z)$, so that there exists
$\beta\in{\mathbb C}$ with $h_1(z)=\alpha h_2(z)+\beta$.

Finally, substituting in equation~(\ref{eq:ThmEquiEnd2}), one has
that $b_1(z)=b_2(z)$.

Let us now prove the converse. Using the formula~(\ref{eq:rhoLk-hbc}), it is  straightforwardly checked that
    \begin{gather*}
    \rho_2(L_{-1})   =   \frac1{\alpha}\rho_1(L_{-1}),
    \\
    \rho_2(L_{i}) =
    h_2(z)^{i} \left(\frac{h_2(z)}{\alpha}\rho_1(L_{-1})-(i+1)c\right)\qquad\forall\,  i\geq -1.
    \end{gather*}
Our f\/irst task is to prove that $\rho_2(L_{i})$ is a linear combination of $\{\rho_1(L_k)\}$. Let us compute $\rho_2( L_i )$
\begin{gather*}
    \rho_2( L_i )  =
    (\alpha h_1(z)+\beta)^{i}\left(\frac{\alpha h_1(z)+\beta}{\alpha}\rho_1(L_{-1}) - (i+1)c\right)
  \\
\hphantom{\rho_2( L_i )}{} =
    \alpha^i\left(h_1(z)+\frac{\beta}{\alpha}\right)^i\left(h_1(z)\rho_1(L_{-1})+\frac{\beta}{\alpha}\rho_1(L_{-1}) -(i+1)c\right)
  \\
\hphantom{\rho_2( L_i )}{} =
    \alpha^i\left(\sum_{j=0}^i\binom{i}{j}h_1(z)^j(\frac{\beta}{\alpha})^{i-j}\right)\left(h_1(z)\rho_1(L_{-1})+\frac{\beta}{\alpha}\rho_1(L_{-1}) -(j+1)c+(j-i)c\right)
  \\
\hphantom{\rho_2( L_i )}{} =
    \alpha^i \sum_{j=0}^i\binom{i}{j}\left(\frac{\beta}{\alpha}\right)^{i-j} \left(\rho_1(L_j)+h_1(z)^j\left(\frac{\beta}{\alpha}\rho_1(L_{-1}) (j-i)c\right)\right)
  \\
\hphantom{\rho_2( L_i )}{}=
    \alpha^i \sum_{j=0}^i\binom{i}{j}\left(\frac{\beta}{\alpha}\right)^{i-j} \left(\rho_1(L_j)+\frac{\beta}{\alpha}\big(\rho_1(L_{j-1})+j c h_1(z)^{j-1}\big)+ (j-i)c h_1(z)^j\right).
    \end{gather*}
Bearing in mind that the term
\begin{gather*}
    \sum_{j=0}^i\binom{i}{j}\left(\frac{\beta}{\alpha}\right)^{i-j}\left(\frac{\beta}{\alpha}\big( j c h_1(z)^{j-1}\big)+ (j-i)c h_1(z)^j\right)
\end{gather*}
vanishes identically, the above expression yields
\begin{gather*}
    \rho_2( L_i )   =
    \alpha^i \sum_{j=0}^i\binom{i}{j}\left(\frac{\beta}{\alpha}\right)^{i-j} \left(\rho_1(L_j)+\frac{\beta}{\alpha}\rho_1(L_{j-1})\right)
   \\
   \hphantom{\rho_2( L_i )}{} =
    \left(\frac{\beta}{\alpha}\right)^{i+1} \rho_1(L_{-1}) +
    \sum_{j=0}^{i-1}\left(\binom{i}{j}+ \binom{i}{j+1}\right)\left(\frac{\beta}{\alpha}\right)^{i-j}\rho_1(L_j) + \rho_1(L_i).
    \end{gather*}
This explicit expression shows at once that $\operatorname{Im}\rho_2=\operatorname{Im}\rho_1$.

Finally, consider
 \begin{gather*}%\label{eq:phialphabeta}
    \phi(L_i )  :=  \left(\frac{\beta}{\alpha}\right)^{i+1} L_{-1} +
    \sum_{j=0}^{i-1}\binom{i+1}{j+1}\left(\frac{\beta}{\alpha}\right)^{i-j} L_j  +  L_i.
    \end{gather*}
The fact that $\phi=(\rho_1\vert_{\operatorname{Im}\rho_1})^{-1}\circ\rho_2$ implies that
$\phi$ is a Lie algebra automorphism of ${\mathcal W}^+$ and
that $\rho_2=\rho_1\circ \phi$.
\end{proof}

\begin{Example}[$m$-reduced KP hierarchy]\label{ex:VirSubalg}
In particular, it is easy to check that given an action $(V,\rho)$  and a non-zero integer number $m$, the map $\rho_m(L_k):= \frac{1}{m}\rho(L_{mk})$ gives rise to another action of ${\mathcal W}^+$ on $V$
\begin{gather*}
    [\rho_m(L_i), \rho_m(L_j)]  =
    \frac{1}{m^2}[\rho(L_{mi}), \rho(L_{mj})]
    =
    \frac{1}{m^2}(mi-mj)\rho(L_{m(i+j)})
    \\
    \hphantom{[\rho_m(L_i), \rho_m(L_j)]}{} =
    (i-j)\frac{1}{m}\rho(L_{m(i+j)})
    =
    (i-j)\rho_m(L_{i+j}).
    \end{gather*}
Moreover, an easy computation yields
\begin{gather*}
    \rho_m(L_k) =
    \frac1{m}\rho(L_{m k})  =
    -\frac{h(z)^{m k+1}}{m h'(z)} \partial_z - \frac1{m}  (m k +1)c h(z)^{m k} + \frac{h(z)^{m k +1}}{m h'(z)} b(z).
\end{gather*}
Thus, this action corresponds to the data $\tilde h(z)=h(z)^m$, $\tilde c=c$ and $\tilde b(z)=(1-m) c \frac{h'(z)}{h(z)} + b(z)$.
\end{Example}

\subsection{Conjugation}\label{subsec:conjugation}

Let us now introduce a generalization of the notion of conjugated
action that will be needed later on.

Let us begin with an example. Assume we have a dif\/ferential equation $P\psi(z)=0$ and we want to solve it by virtue of a replacement $\psi(z)=v(z)\phi(z)$, for a given function $v(z)$; that is, we shall solve $P(v(z)\phi(z))=0$ for an unknown function $\phi(z)$. This is equivalent to solving $(v(z)^{-1}\circ P\circ v(z))\phi(z)=0$, where $v(z)$ is regarded as an operator; namely, the homothety of ratio $v(z)$. For instance, if $P$ is a f\/irst-order dif\/ferential operator with symbol $\sigma(P)$, it holds that
\begin{gather*}
    v(z)^{-1} ( P \psi (z) ) =
    \big(v(z)^{-1}\circ P \circ v(z)\big)\phi(z)   =
    \left(P +\sigma(P) \frac{v'(z)}{v(z)}\right)\phi(z).
\end{gather*}
Therefore, solving the dif\/ferential equation $P\psi(z)=0$ is equivalent to solving
\begin{gather*}
\left(P +\sigma(P) \frac{v'(z)}{v(z)}\right)\phi(z)=0.
\end{gather*}
Observe that, from the point of view of analysis, the space of functions in which~$\psi$ and~$\phi$ lie might be dif\/ferent (e.g., dif\/ferent convergence domains, etc.). In particular, conjugation is an instance of gauge transformation.

Let us recall from \cite{HGPM-sgl} the def\/inition of the group of
semilinear transformations and some of its properties. The group
of semilinear transformations of a f\/inite-dimensional
${\mathbb C}((z))$-vector space~$V$, denoted by~$\operatorname{SGl}_{{\mathbb C}((z))}(V)$,
consists of ${\mathbb C}$-linear automorphisms $\gamma\colon V\to V$ such that
there exists a ${\mathbb C}$-algebra automorphism of~${\mathbb C}((z))$, $g$,
satisfying
    \begin{gather}\label{eq:semilinearGammaG}
    \gamma(f(z)\cdot v) =  g(f(z))\cdot \gamma(v)\qquad \forall\,  f(z)\in{\mathbb C}((z)),  \quad  v\in V,
    \end{gather}
and, therefore, $\operatorname{SGl}({\mathbb C}((z)))=\operatorname{Aut}_{{\mathbb C}\text{-alg}}{\mathbb C}((z)) \ltimes
{\mathbb C}((z))^*$.

The Lie algebra of $\operatorname{SGl}_{{\mathbb C}((z))}(V)$ consists of f\/irst-order
dif\/ferential operators on $V$ with scalar symbol,
${\mathcal D}^1_{{\mathbb C}((z))/{\mathbb C}}(V,V)$, and the symbol coincides with the map
induced by the group homomorphism that sends $\gamma$ to $g$
(related by equation~(\ref{eq:semilinearGammaG})) between their
Lie algebras.

\begin{Theorem}\label{thm:HomModAutxC((z))}
The space $\operatorname{Hom}_{\text{\rm Lie-alg}}({\mathcal W}^+,{\mathcal D}^1)\setminus
\{0\}$ carries an action of the group
$\operatorname{SGl}({\mathbb C}((z)))$ by conjugation and the quotient space is ${\mathbb Z}\times
{\mathbb C}\times \big({\mathbb C}((z))/{\mathbb Z} z^{-1}+{\mathbb C}[[z]]\big)$.
\end{Theorem}

\begin{proof}
Let us begin studying the action of the automorphism group
$G:=\operatorname{Aut}_{{\mathbb C}\text{-alg}}{\mathbb C}((z))$ (for a study and
applications of this group, see~\cite{MP2}). Let us denote elements
of~$G$ with big  Greek letters ($\Phi,\Psi,\ldots)$ and, for each
of them, let the corresponding small Greek letter denote the image
of~$z$; that is
\begin{gather*}
    \Phi(f(z)) =   f(\phi(z))
\end{gather*}
and observe that ${\mathfrak v}(\phi(z))=1$ in order for $\Phi$ to be an isomorphism.

We consider the action of $G$ on the space of actions by conjugation; i.e.,
\begin{gather*}
    (\Phi,\rho)  \mapsto   \rho^\Phi \qquad\text{where} \quad \rho^\Phi(L_k):= \Phi\circ\rho(L_k)\circ\Phi^{-1}\quad \forall\,  k.
\end{gather*}
Let us check that this def\/inition makes sense. Let $\rho$ be given by a triple $(h(z),c,b(z))$. It is straightforward that
    \begin{gather*}
    \rho^\Phi(L_{-1})f(z)  =
    \big(\Phi\circ\rho(L_{-1})\circ\Phi^{-1}\big)f(z) =
    \Phi\left(\left(-\frac{1}{h'(z)}\partial_z+\frac{b(z)}{h'(z)}\right)f(\phi^{-1}(z))\right)  \\
    \hphantom{\rho^\Phi(L_{-1})f(z)}{}
    =
    \left(-\frac{(\phi^{-1})'(\phi(z))}{h'(\phi(z))}\partial_z +
    \frac{b(\phi(z))}{h'(\phi(z))}\right)f(z).
    \end{gather*}
Note that expanding and derivating the identity $\Phi\circ\Phi^{-1}(z)=z$, one gets that $(\phi^{-1})'(\phi(z))\cdot \phi'(z)=1$ and, thus
\begin{gather*}
    \frac{(\phi^{-1})'(\phi(z))}{h'(\phi(z))}   =
    \frac{1}{\partial_z h(\phi(z))}.
\end{gather*}
Summing up, the transformation $\Phi$ acts on triples as follows
\begin{gather*}
 (\Phi , (h(z),c,b(z)) )  \mapsto   (h(\phi(z)),c,b(\phi(z))).
\end{gather*}

Second, we study the action of ${\mathbb C}((z))^*$. Bearing in mind the
discussion of the beginning of this subsection, we consider the
action
\begin{gather*}
    (s(z), \rho ) \mapsto   \big( s(z) \circ \rho \circ   s(z)^{-1}\big),
\end{gather*}
so that, in terms of triples, it holds that
\begin{gather*}
 (s(z), (h(z),c,b(z))  ) \mapsto   \left(h(z),c,b(z)-\frac{s'(z)}{s(z)}\right).
\end{gather*}

One checks easily that the f\/irst def\/ined action intertwines the
second one; that they yield an action of the semidirect product $G
\ltimes {\mathbb C}((z))^*$; and that the quotient space under this action
is ${\mathbb Z}\times {\mathbb C}\times ({\mathbb C}((z))/{\mathbb Z} z^{-1}+{\mathbb C}[[z]])$, where an
action $\rho$, corresponding to a triple $(h(z),c,b(z))$ is mapped
to $({\mathfrak v}(h(z)), c, \bar b(z))$ ($\bar b(z)$ being the
equivalence class of~$b(z)$).
\end{proof}

Nevertheless, the conjugation also makes sense if $v(z)$ is
replaced by any  linear
operator on the space of functions such that $\frac{v'(z)}{v(z)}$
can be identif\/ied with an element in ${\mathbb C}((z))$. This is the case of
the example at the beginning of this subsection, but it also holds for functions~$v(z)$
admitting an asymptotic expansion at $0$. For instance, for the formal
expression $v(z):= \exp(\int s(z)\operatorname{d}z)$ where
$s(z)\in{\mathbb C}((z))$, the quotient $\frac{v'(z)}{v(z)}$ will be
identif\/ied with $s(z)$. Indeed, the conjugation by
$\exp(-\frac23z^{-3})$ was used in~\cite{KacSchwarz} when solving
a dif\/ferential equation. For another example, let us consider
$v(z)$ to be a solution of the second-order dif\/ferential equation
$v''(z)+\frac{1}{2}S(h)v(z)=0$, where $S(h)$ denotes the
Schwarzian derivative of $h$, such that $\frac{v'(z)}{v(z)}\in
{\mathbb C}((z))$, which holds true in many cases (e.g., when
$S(h)\in{\mathbb C}((z))$).

It is worth
noticing that once $\frac{v'(z)}{v(z)}$ is thought of as an
element of ${\mathbb C}((z))$, $\frac{v''(z)}{v(z)}$ will be identif\/ied
with $\big(\frac{v'(z)}{v(z)}\big)^2 +
\big(\frac{v'(z)}{v(z)}\big)'\in{\mathbb C}((z))$. By abuse of notation, we
def\/ine $\operatorname{d}\log v(z):=
\frac{v'(z)}{v(z)}$.

Thus, for $P\in{\mathcal D}^1({\mathbb C}((z)))$ and $v(z)$ as above, we
consider another f\/irst-order dif\/ferential operator $P^v\in{\mathcal
D}^1({\mathbb C}((z))\otimes_{{\mathbb C}}{\mathbb C} v(z))$ def\/ined by
\begin{gather*}
    P^v (f(z)\otimes v(z))  :=
    \left(\left(P +\sigma(P) \frac{v'(z)}{v(z)}\right)(f)\right)\otimes v(z).
\end{gather*}
The induced map from ${\mathcal D}^1({\mathbb C}((z)))$ to ${\mathcal
D}^1({\mathbb C}((z))\otimes_{{\mathbb C}}{\mathbb C} v(z))$ is a Lie algebra homomorphism.

\begin{Definition}
The conjugated action of $(V,\rho)$ by $v(z)$ is the pair
$(V^v,\rho^v)$, consisting of the $1$-dimensional ${\mathbb C}((z))$-vector
space $V^v:=V\otimes_{{\mathbb C}}{\mathbb C} v(z)$ together with the action def\/ined
by
    \begin{gather}\label{eq:conjugationRhoVequal}
    \rho^v(L_k) (f(z)\otimes v(z))  :=
    \left(\rho(L_k) (f(z))    +  \sigma(\rho(L_k)) f(z) \frac{v'(z)}{v(z)} \right)\otimes v(z).
    \end{gather}
\end{Definition}

In particular, if the data $h(z)$, $c$, $b(z)$ def\/ine an action~$\rho$,
then $h(z)$, $c$, $b(z)-\frac{v'(z)}{v(z)}$ def\/ine~$ \rho^v$.

\begin{Remark}
Our generalized notion of conjugation allows us to consider other actions of~${\mathbb C}((z))$. For instance, if we consider $s(z)\in{\mathbb C}((z))$ acting by conjugation by $\exp(s(z))$ (resp.\  $\exp(\int s(z))$), then the quotient space by $G\ltimes {\mathbb C}((z))$ is ${\mathbb Z}\times {\mathbb C}\times {\mathbb C}/{\mathbb Z}$ (resp.\ ${\mathbb Z}\times {\mathbb C}$, which is related to Theorem~\ref{thm:IsomClassVirasoro}).
\end{Remark}

\subsection{Central extensions}\label{subsec:CentralExt}

Recall from~\cite{AMP} that, given a pair $(V,V_+)$ consisting of a ${\mathbb C}$-vector space and a subspace of it, there exists a ${\mathbb C}$-scheme whose rational points correspond to the vector subspaces $W\subseteq V$ such that
\begin{gather*}
W\cap V_+ \qquad \mbox{and}\qquad  V/W+V_+
\end{gather*}
are f\/inite-dimensional vector spaces over ${\mathbb C}$. It is called the Sato Grassmannian or inf\/inite Grassmannian of~$V$ (see~\cite{SegalWilson} for an analytical approach).

Further, the group $\operatorname{SGl}({\mathbb C}((z)))$ acts on the Sato Grassmannian of $V$,
$\operatorname{Gr}(V)$, preserving the determinant bundle $\operatorname{Det}_V$. This has two
signif\/icant consequences. Firstly,  this group also acts on the
projectivization of $H^0(\operatorname{Gr}(V),\operatorname{Det}_V^*)$, making the Pl\"{u}cker
embedding equivariant
\begin{gather*}
    \operatorname{Gr}(V)  \hookrightarrow  {\mathbb P} H^0\big(\operatorname{Gr}(V),\operatorname{Det}_V^*\big)^*.
\end{gather*}
Second, the group of
automorphisms of $\operatorname{Det}_V$ that lift automorphisms of~$\operatorname{Gr}(V)$
def\/ined by elements of $\operatorname{SGl}_{{\mathbb C}((z))}(V)$; i.e.,
\begin{gather*}
    \widetilde{\operatorname{SGl}}_{{\mathbb C}((z))}(V)
    := \left\{
    \raisebox{18pt}{{ \xymatrix@R=13pt{
        \operatorname{Det}_V^* \ar[r]^{\sim} \ar[d] & \operatorname{Det}_V^*  \ar[d]
        \\
        \operatorname{Gr}(V) \ar[r]_{\sim}^g & \operatorname{Gr}(V)}}}
        \quad\text{where} \ g\in \operatorname{SGl}_{{\mathbb C}((z))}(V)
    \right\}
\end{gather*}
def\/ines a canonical central extension
    \begin{gather}\label{eq:CentralExtSGl}
    1\to {\mathbb C}^*{\mathbf c} \to \widetilde{\operatorname{SGl}}_{{\mathbb C}((z))}(V) \to \operatorname{SGl}_{{\mathbb C}((z))}(V) \to 1,
    \end{gather}
where ${\mathbf c}$ is the central element (\emph{central
charge}).

Now, let us assume that an action $(V,\rho)$ of ${\mathcal W}^+$
is given. Since $\rho$ can be extended to ${\mathcal W}$, we may
consider the pullback of the central extension of Lie algebras induced by
(\ref{eq:CentralExtSGl}) and obtain a central extension
$\widetilde {\mathcal W}$ of ${\mathcal W}$
    \begin{gather*}%\label{eq:CentralExtW}
    0\to {\mathbb C}{\mathbf c} \to \widetilde {\mathcal W} \to {\mathcal W} \to 0,
    \end{gather*}
where we recall that $\widetilde {\mathcal W} = {\mathcal W}\oplus {\mathbb C}{\mathbf c}$, as vector spaces. Its Lie bracket will be computed below.

\begin{Remark}\label{rem:KPflowSGLflows-gauge}
Recall that the multiplicative group ${\mathbb C}((z))^*$ acts by homotheties on $\operatorname{Gr}(V)$ and by conjugation (i.e., gauge transformation) on $\operatorname{SGl}_{{\mathbb C}((z))}(V)$. Furthermore, the action map
\begin{gather*}
    \operatorname{SGl}_{{\mathbb C}((z))}(V)\times\operatorname{Gr}(V)  \longrightarrow  \operatorname{Gr}(V)
\end{gather*}
is equivariant with respect to these actions. Thus, for f\/ixed
$U\in\operatorname{Gr}(V)$, the map induced by the orbit map at the level of
tangent spaces
\begin{gather*}
    T_{\operatorname{Id}}\operatorname{SGl}_{{\mathbb C}((z))}(V)  \simeq
     {\mathcal D}^1_{{\mathbb C}((z))/{\mathbb C}}(V,V)   \longrightarrow  T_U\operatorname{Gr}(V) \simeq
\operatorname{Hom}(U,V/U)
\end{gather*}
intertwines gauge transformations (on the l.h.s.) and KP f\/lows (on
the r.h.s.). In fact,  requirements~(ii$'$) and~(iii$'$) of the
introduction, $z^{-2}U\subset U$ and $\tilde L_kU\subseteq U$,
should be understood at the tangent space~$T_U\operatorname{Gr}( V)$.
\end{Remark}

In fact, once  $(V,\rho)$  is given, the above discussion allows us to obtain a natural map at the level of tangent spaces
    \begin{gather}\label{eq:WActsTangentGr}
    \begin{split} &
    \xymatrix@R=0pt{
    {\mathcal W}^+ \ar[r]^-{\rho} & {\mathcal D}^1_{{\mathbb C}((z))/{\mathbb C}}(V,V) \ar[r] & \operatorname{Hom}(U,V/U),
    \\
    L  \ar@{|->}[rr] & & U\hookrightarrow V\overset{\rho(L)}\to V \to V/U,
    }\end{split}
    \end{gather}
where $U$ is a point in $\operatorname{Gr}(V)$ and we observe that this map can
be extended to ${\mathcal W}$ and that the equivariance of the Pl\"{u}cker morphism yields a natural map
\begin{gather*}
    \xymatrix@R=0pt{
    \widetilde {\mathcal W}  \ar[r]^-{\hat\rho} &
    H^0(\operatorname{Gr}(V),\operatorname{Det}_V^*)^* \ar[r]^-{\sim} &
    {\mathbb C}[[t_1,t_2,\ldots]],
    \\
    L \ar@{|->}[r] &
    { \text{evaluation at }\rho(L)(U)} \ar@{|->}[r] &
    \hat\rho(L)\tau_U(t)}
\end{gather*}
($\mathbf c$, being the central element, acts as a homothety).

\begin{Remark}
It is worth pointing out that $H^0(\operatorname{Gr}(V),\operatorname{Det}_V^*)$, also written
as $\Lambda^{\frac{\infty}{2}}V$ in the literature, corresponds to
the formulation in  terms of  fermions and that the isomorphism on
the r.h.s.\ is the bosonization morphism (see, e.g.,~\cite{KacRa}).
Observe that the rightmost term of the previous map,
${\mathbb C}[[t_1,t_2,\ldots]]$, does not depend on $V$ since
the $\tau$-function of a point $U$  depends only on the
coordinates of vectors of $U$ w.r.t.\ a basis of $V$ but not on the
basis (e.g., \cite[Lemma~4.2.]{Kon}).
\end{Remark}

Our next task is concerned with the explicit computation of the
Lie bracket of $\widetilde {\mathcal W}$ or, tantamount to this,
with the classif\/ication of $\widetilde {\mathcal W}$ as a central
extension.
Since $H^2({\mathcal W},{\mathbb C})\simeq {\mathbb C}$ and it is generated  by the
cocycle of the Virasoro algebra,
$\operatorname{Vir}(L_r,L_s)=\delta_{r,-s}\frac{r^3-r}{12}$, it
follows that the $2$-cocycle classifying the central extension
$\widetilde {\mathcal W}$ has to be a multiple of~$\operatorname{Vir}$.

It is known that the cocycle of the extension~(\ref{eq:CentralExtSGl}) is the so-called Gelfand--Fuchs cocycle whose general expression reads as follows (see, for instance,~\cite{KacRa,KacPe})
\begin{gather}\label{eq:2cocycle}
    \Psi\big(a(z)\partial_z^r, b(z)\partial_z^s\big) :=
    \frac{r! s!}{(r+s+1)!} \mathop{\operatorname{Res}}\limits_{z=0}
    \big(\partial_z^{s+1}a(z)
      \partial_z^{r}b(z)\big)\operatorname{d}z.
    \end{gather}
This means that $\widetilde {\mathcal W}$, generated by
$\{L_k,{\text{\textbf c}}\}$, carries the Lie bracket
\begin{gather*}
    [L_r,L_s]:=(r-s)L_{r+s}+\Psi(L_r,L_s){\text{\textbf c}}  .
\end{gather*}

\begin{Theorem}\label{thm:IsomClassVirasoro}
Let $(V,\rho)$ be given.
The isomorphism class of the $2$-cocyle of $\widetilde {\mathcal
W}$, $\Psi$, is
\begin{gather*}
    - 2\big(1-6c+6c^2\big) {\mathfrak v}(h(z))  \operatorname{Vir},
\end{gather*}
where $h(z)$, $c$ are given by Theorem~{\rm \ref{thm:VirasoroSubAlg}}.
\end{Theorem}

\begin{proof}
Let us recall what a $2$-coboundary looks like. A $2$-cocycle
\begin{gather*}
\alpha\colon \ {\mathcal W}\times{\mathcal W}\to {\mathbb C}
\end{gather*}
is a $2$-coboundary if\/f there exists a linear map $f\colon {\mathcal W}\to {\mathbb C}$ such that
\begin{gather*}
    \alpha(L,L')   =  f([L,L']).
\end{gather*}
Thus, in our case, since $\{L_r\,\vert\, r\in {\mathbb Z}\}$ is a basis of
${\mathcal W}$, a function $f\colon {\mathcal W}\to {\mathbb C}$ is expressed as
$f=\sum_r f_r\Delta_r$, where $f_r\in{\mathbb C}$ and $\{\Delta_r\}$ is the
dual basis; i.e., $\Delta_r(L_s)=\delta_{r,s}$ (and, hence,
$f_r=f(L_r)$). In the case of ${\mathcal W}$, it follows that~$\alpha$ is a $2$-coboundary if\/f there exist constants~$\{f_r\}$
such that $\alpha(L_r,L_s) = (r-s)\cdot f_{r+s}$; or,
equivalently, there exists a function $\bar f\colon {\mathbb Z}\to {\mathbb C}$ such that
\begin{gather*}
    \alpha(L_r,L_s)  =   (r-s)  \bar f(r+s).
\end{gather*}
Observe that the isomorphism class of~$\Psi$ is a multiple of
$\operatorname{Vir}$ and it suf\/f\/ices  to compute the quotient of
$\Psi(\rho(L_{r}),\rho(L_{-r}))$ by~$\frac{r^3-r}{12}$, where
$\Psi$ is given by equation~(\ref{eq:2cocycle}). We do not care
about the remainder, since it will take the form of a
$2$-coboundary.

Note that the isomorphism class of~$\Psi$ does not depend on the
conjugacy class of~$\rho$, so that it can be assumed that
$V={\mathbb C}((z))$ and, further, $h(z)$ and $c$ are invariant under
conjugation (see Section~\ref{subsec:conjugation}). Thus,
given $h(z)$, $c$, $b(z)$~and $\rho$ in the form of
Theorem~\ref{thm:VirasoroSubAlg}, one has that
$\Psi(\rho(L_{r}),\rho(L_{-r}))$ is the residue at $z=0$ of the
following expression
    \begin{gather*}
    \frac{1}{3!}
    \left(\partial_z^{2}\left(\frac{-h(z)^{r+1}}{h'(z)}\right)
      \partial_z\left(\frac{-h(z)^{-r+1}}{h'(z)}\right)\right)
\\
\qquad{} +
    \frac{1}{2!}
    \left(\partial_z\left(\frac{-h(z)^{r+1}}{h'(z)}\right)
      \partial_z\left((1-r)c  h(z)^{-r}+ \frac{h(z)^{-r+1}}{h'(z)}b(z)\right)\right)
  \\
  \qquad{} +
    \frac{1}{2!}
    \left(\partial_z^2\left(-(r+1)c h(z)^{r}+\frac{h(z)^{r+1}}{ h'(z)}b(z)\right)
      \left(\frac{-h(z)^{-r+1}}{h'(z)}\right) \right)
  \\
  \qquad {} +
    \left(\partial_z\left(-(r+1)c  h(z)^{r}+\frac{h(z)^{r+1}}{ h'(z)}b(z)\right)
      \left((1-r)c  h(z)^{-r}+\frac{h(z)^{-r+1}}{h'(z)}b(z)\right) \right).
    \end{gather*}
Expanding this expression and dividing by $\frac{r^3-r}{12}$, one obtains the following quotient
\begin{gather*}
    -2\big(1-6c+6c^2\big)\frac{h'(z)}{h(z)},
\end{gather*}
whose residue is equal to
\begin{gather*}
    -2\big(1-6c+6c^2\big) {\mathfrak v}(h(z))
\end{gather*}
and the result follows.
\end{proof}

\begin{Remark}\label{rem:densities}
Note that when~$c$ is an integer, the central extension admits an
alternative geometric construction, namely, the one induced by the
action of~${\mathbb C}((z))$ on~$\operatorname{Gr}({\mathbb C}((z))(\operatorname{d}z)^{\otimes c})$ (i.e., on the
space of $c$-densities). This approach was used in~\cite{MP2} to
of\/fer a geometric formalism of the bosonic string. It is worth
pointing out that the actions attached to the data~$(h(z),c,b(z))$
and to $(h(z),1-c,b(z))$ are isomorphic.
\end{Remark}

\begin{Remark}
Let $\widehat{\mathcal D}$ (a.k.a.~${\mathcal W}_{1+\infty}$)
denote the unique non-trivial central extension of the Lie algebra
of dif\/ferential operators on the circle. Owing to its presence in
several topics, the authors of~\cite{FKRW} carried out an in-depth
study of its representations with the help of the theory of vertex
operator algebras (voa). In this context, they consider two
$1$-parameter families of Virasoro algebras (see
\cite[equation~(1.7)]{FKRW}) $\{L^+_k(\beta)\,\vert\, k\in{\mathbb Z}\}$ and
$\{L^-_k(\beta)\vert k\in{\mathbb Z}\}$, where
    \begin{gather*}
    \big\{ L^+_k(\beta) = -z^{k+1}\partial_z - \beta(k+1) z^k
    \, \vert\, k\geq -1\big\},
    \\
    \big\{ L^-_k(\beta) = -z^{k+1}\partial_z - (k+\beta(-k+1)) z^k
    \, \vert\, k\geq -1\big\},
    \end{gather*}
which eventually play a relevant role in their results on
conformal voa (see \cite[Theorem~3.1.]{FKRW}). Let us check that,
for each $\beta$, these Virasoro algebras f\/it in our approach. For
the case $\{L^+_k(\beta)\}$ the system is
    \begin{gather*}
 -\frac{h(z)^{i+1}}{h'(z)} = - z^{i+1},
\qquad
     - (i+1) c  h(z)^i + \frac{b(z)}{h'(z)} h(z)^{i+1} =  - \beta(i+1)
    z^i
    \end{gather*}
and it yields $h(z)=z$, $c=\beta$ and $b(z)=0$. They show that
the central charge, ${\mathbf c}(\beta)$, satisf\/ies ${\mathbf
c}(\beta)=-2(6\beta^2-6\beta+1) {\mathbf c}(1)$, and indeed this
agrees with Theorem~\ref{thm:IsomClassVirasoro}. For the second
family, $\{L^-_k(\beta)\}$, one obtains $h(z)=z$, $c=1-\beta$ and
$b(z)=(1-2\beta)z^{-1}$.
\end{Remark}

%%%%%%%%%%%%%%%%%%%%%%%%%%%%%%%%%%%%%%%%%%%%%
\section{KdV and the string equation}\label{sec:SolutionsTAUVir&KP}
%%%%%%%%%%%%%%%%%%%%%%%%%%%%%%%%%%%%%%%%%%%%%

%%%%%%%%%%%%%%%%%%%%%%
\subsection{Stabilizer}\label{subsec:stab}

Following the notation of~Section~\ref{subsec:ComputationVirasoroSubalg}, let $(V,\rho)$ be an action of ${\mathcal W}^+$.

Let $U\subset V$ denote a ${\mathbb C}$-vector subspace and let $A_U$ denote its stabilizer; that is
\begin{gather*}
    A_U :=  {\operatorname{Stab}}(U) =  \{f\in {\mathbb C}((z)) \,\vert\, fU\subseteq U\}.
\end{gather*}

We say that $U$ is $L_{-1}$-\emph{stable} or \emph{stable under the action of}
$L_{-1}$ when $\rho(L_{-1})U\subseteq U$. Similarly, we say that
$U$ is ${\mathcal W}^+$-stable when $\rho(L)U\subseteq U$ for all
$L\in {\mathcal W}^+$; or, what is tantamount to this, $\rho(L_k)
U\subseteq U$ for all  $k\geq -1$.

Let us f\/ix the following notation. Let $V_+\subseteq V$ denote a ${\mathbb C}[[z]]$-submodule of $V$ and, as above, let $\mathfrak v$ be the valuation def\/ined by~$z$.

\begin{Theorem}\label{thm:stab=C[h]orC}
Let $(V,\rho)$ be an action of ${\mathcal W}^+$ and let $h(z)$ be given as in Theorem~{\rm \ref{thm:VirasoroSubAlg}}. Let~$U$ be a subspace of~$V$ such that $U\neq (0)$, $U\cap V_+$ is finite-dimensional and  ${\mathfrak v}(h(z)) < 0$.

If $U$ is $L_{-1}$-stable and $A_U\neq {\mathbb C}$, then $U$ is ${\mathcal W}^+$-stable  and  $A_U={\mathbb C}[h(z)]$.
\end{Theorem}

\begin{proof}
First, let us show that ${\mathfrak v}(f(z))<0$ for each $f(z)\in
A_U$ non-constant. Indeed, if ${\mathfrak v}(f(z))>0$ then $U\cap
V_+$ cannot be f\/inite-dimensional since $U\neq (0)$. On the other
hand, if ${\mathfrak v}(f(z))=0$, then $\bar f(z):=f(z)-f(0)$
belongs to $A_U$ and ${\mathfrak v}(\bar f(z))>0$, which again
contradicts the hypotheses. Thus, it must hold that ${\mathfrak
v}(f(z))<0$.

We shall now prove that $A_U\subseteq {\mathbb C}[h(z)]$. From the previous
paragraph, let us take $f(z)\in A_U\setminus {\mathbb C}[h(z)]$ such that
${\mathfrak v}(f(z))$ attains the value $\max\{{\mathfrak v}(f(z))
\,\vert \, f(z)\in A_U\setminus {\mathbb C}[h(z)]\}$. Since $U$ is stable under
$L_{-1}$ and under the multiplication by $f(z)$, it follows that
$[f(z),\rho(L_{-1})]=\frac{f'(z)}{h'(z)}\in A_U$. Note that
${\mathfrak v}(\frac{f'(z)}{h'(z)})={\mathfrak v}(f(z))-{\mathfrak
v}(h(z)) >{\mathfrak v}(f(z))$, since ${\mathfrak v}(h(z))$ is negative and ${\mathfrak v}(f(z))\neq 0$. Bearing in mind  that $f(z)$
is such that ${\mathfrak v}(f(z))$ is maximal among those elements
of $A_U\setminus {\mathbb C}[h(z)]$, we have that
\begin{gather*}
    \frac{f'(z)}{h'(z)} \in {\mathbb C}[h(z)]
\end{gather*}
and, thus, $f(z)\in {\mathbb C}[h(z)]$. That is, $A_U \subseteq{\mathbb C}[h(z)]$.

Let us see that $A_U={\mathbb C}[h(z)]$. Since $A_U\neq{\mathbb C}$, let $p(x)$ be a
non-constant polynomial of minimal degree such that $p(h(z))\in
A_U$. Similar to the above, one has that
$[p(h(z)),\rho(L_{-1})]=p'(h(z))\in A_U$ and, thus, $p'(x)$ must
be constant and $p(x)$ is of the form $a x + b$. Therefore,
${\mathbb C}[h(z)]={\mathbb C}[p(h(z))]\subseteq A_U\subseteq {\mathbb C}[h(z)]$.

It remains to show that, in the case $A_U={\mathbb C}[h(z)]$, $U$ is ${\mathcal W}^+$-stable. This follows from the fact that $h(z)U\subseteq U$ and from the relation  $\rho(L_i) = h(z)^i(h(z)\rho( L_{-1})-(i+1)c)$ for all $i\geq 0$.
\end{proof}

\begin{Remark}
The hypothesis ${\mathfrak v}(h(z)) < 0$ is satisf\/ied in most relevant cases. For instance, if $U$ is ${\mathcal W}^+$-stable and  $U\cap V_+$ is f\/inite-dimensional, then ${\mathfrak v}(h(z)) \leq 0$. Indeed, if ${\mathfrak v}(h(z)) > 0$, then the operator  $ h(z)^i(h(z)\rho( L_{-1})-(i+1)c)$ would raise the order for all $i\gg 0$ and, since $U\neq (0)$, $U\cap V_+$ could not be of f\/inite dimension.
\end{Remark}

\begin{Remark}
This result clearly unveils  why the string equation together with
the KdV hierarchy (i.e., $L_{-1}$-stability and $h(z)=z^{-2}$)
imply the Virasoro constraints (i.e., ${\mathcal W}^+$-stability).
On the other hand, the existence of points in the Grassmannian
with trivial stabilizer that are $L_{-1}$-stable but not
${\mathcal W}^+$-stable is known \cite[Section~8]{Mulase-AlgKP}.
\end{Remark}

Given an action $(V,\rho)$, let us introduce the \emph{first-order
stabilizer} of a subspace $U\subset V$ as
    \begin{gather}\label{eq:A1U}
    A_U^1  :=  \big\{D\in {\mathcal D}^1_{{\mathbb C}((z))/{\mathbb C}}(V,V) \,\vert \, D(U)\subseteq U\big\}.
    \end{gather}
Then, for a ${\mathcal W}^+$-stable subspace $U\subset V$ such that $A_U={\mathbb C}[h(z)]$, there is a canonical exact sequence of Lie algebras
    \begin{gather}\label{eq:A_U^1split}
    \xymatrix{
    0 \ar[r] &A_U \ar[r] &
    A_U^1  \ar[r]  &
    \operatorname{Der}_{{\mathbb C}}(A_U) \ar[r] &
    0,
    }
    \end{gather}
and $\rho$ induces a splitting of it. To prove this, one uses the
sequence
\begin{gather*}
      0\to \operatorname{End}_{{\mathbb C}((z))}V \to {\mathcal D}^1_{{\mathbb C}((z))/{\mathbb C}}(V,V)\rightarrow \operatorname{Der}_{{\mathbb C}}\big({\mathbb C}((z))\big)\to 0 ,
\end{gather*}
together with \cite[Remark~3.9]{HGPM-sgl} and the explicit expressions of the operators $\rho(L_i)$ computed in Theorem~\ref{thm:VirasoroSubAlg}.

\begin{Remark}
Note that if $U$ is $L$-stable for $L\in {\mathcal
D}^1_{{\mathbb C}((z))/{\mathbb C}}(V,V)$, it follows that $L\in A_U\oplus
\operatorname{Im}(\rho)$. Further, given an action $({\mathbb C}((z)),\rho)$ and a ${\mathcal W}^+$-stable
point $U\in\operatorname{Gr}(V)$, then
$A_U^1={\mathbb C}[h(z)]\otimes_{{\mathbb C}}\langle 1,\rho(L_{-1})\rangle$ as a Lie subalgebra
of ${\mathcal D}^1_{{\mathbb C}((z))/{\mathbb C}}(V,V)$. It is worth mentioning the
paper  \cite[Section~2.1]{AMSvM} where the authors study subspaces of
the Sato Grassmannian, which are stable by the multiplication by a
power of $z$ as well as by the action of a f\/irst-order
dif\/ferential operator. Then,  they investigate the matrix integral
representation of the corresponding $\tau$-function.
\end{Remark}

\begin{Remark}\label{rem:BigDiagramGr}
Recall from Remark~\ref{rem:KPflowSGLflows-gauge} that the groups
${\mathbb C}((z))^*$ and $\operatorname{SGl}(V)$ act canonically on the Sato
Grassmannian. For $U$ a point of  $\operatorname{Gr}(V)$, it is well known that
the stabilizer, $A_U$, coincides with the kernel of the
dif\/ferential of the orbit map, $T_{1}{\mathbb C}((z))^*\to T_U\operatorname{Gr}(V)$.
Similarly, in the case of semilinear automorphisms the orbit map
at the tangent spaces reads as follows
\begin{gather*}
    {\mathcal D}^1_{{\mathbb C}((z))/{\mathbb C}}(V,V)  =
    T_{1} \operatorname{SGl} (V)
    \longrightarrow
    \operatorname{Hom}_{{\mathbb C}}(U,V/U),
\end{gather*}
such that the f\/irst-order stabilizer $A_U^1$ is its kernel.
\end{Remark}

\begin{Proposition}\label{prop:LieAlgLift}
Let $({\mathbb C}((z)),\rho)$ be an action associated with a triple
$(h(z),c,b(z))$. Let $U\subset {\mathbb C}((z))$ be a ${\mathcal
W}^+$-stable subspace such that $A_U={\mathbb C}[h(z)]$. Let  $\widetilde
A_U^1$ be the central extension of $A_U^1$  obtained from the
exact sequence~\eqref{eq:CentralExtSGl}. Let $b(z)=\sum_i b_i
z^i$.

If $h(z)=z^{-n}$, $b_{ni-1}=0$ for all $i\geq 1$ and
$b_{-1}=\frac{n+1}2$, then the exact sequence of Lie algebras
\begin{gather*}
        0\to {\mathbb C}{\mathbf c} \to \widetilde A_U^1 \to A_U^1 \to 0
\end{gather*}
splits.
\end{Proposition}

\begin{proof}
Similarly to the proof of Theorem~\ref{thm:IsomClassVirasoro}, we
shall compute $\Psi$, the Gelfand--Fuchs cocycle
(equation~(\ref{eq:2cocycle})), for a basis of~$A_U^1$. Recall
that a basis is given by~$\rho(L_k)$ for~$k\geq-1$ and by~$h(z)^k$
for~$k\geq 0$.

Recall that
\begin{gather*}
    \Psi(\rho(L_r),\rho(L_s))   =
    \big(1-6c+6c^2\big)n\frac{r^3-r}{6}\delta_{r,-s},
\end{gather*}
which vanishes for $r,s\geq -1$.

Secondly, note that $\Psi(h(z)^r,h(z)^s)=- r n \delta_{r,-s}$ and
that, hence, it vanishes for $r,s\geq 0$.

Finally, a straightforward computation yields
    \begin{gather*}
    \Psi\big(\rho(L_r)+   (r+1)c h(z)^r ,h(z)^s\big)
  =
    \Psi\left(-\frac{h(z)^{r+1}}{h'(z)}\partial_z , h(z)^s\right)   +
    \Psi\left(\frac{h(z)^{r+1}}{h'(z)} b(z) ,h(z)^s\right)
  \\
\qquad   =     \frac12  \mathop{\operatorname{Res}}\limits_{z=0} \left(\partial_z\left(-\frac{h(z)^{r+1}}{h'(z)}\right)   \partial_z( h(z)^s)\right)\operatorname{d}z    +
   \mathop{\operatorname{Res}}\limits_{z=0} \left( \partial_z\left(\frac{h(z)^{r+1}}{h'(z)} b(z)\right)
    h(z)^s\right)\operatorname{d}z
     \\
     \qquad{} =
    \frac12(nr-1) n s \delta_{r,-s} + n s b_{n(r+s)-1},
    \end{gather*}
which vanishes for all $r\geq -1$, $s\geq 0$, because of the hypotheses on the coef\/f\/icients~$b_i$.
\end{proof}

\begin{Example}\label{exam:ConstraintResidueb-1}
Let us consider the case where $h(z)=z^{-2}$ and $b_{-1}\in\frac12{\mathbb Z}\setminus{\mathbb Z}$. We will see that there exists $s(z)\in{\mathbb C}((z))$ such that
the sequence for $\widetilde A_{s(z) U}^1$ is split. First of all,
note that for $\bar s(z)=1+\sum\limits_{i>0} s_i z^{2i} \in {\mathbb C}[[z]]^*$,
it holds that the stabilizer of $\bar s(z) U$ is generated by
$h(z)^r$ for $r\geq 0$ and by $\rho^{\bar s(z)}(L_r)$ for $r\geq
-1$ where, as in Section~\ref{subsec:conjugation}, $\rho^{\bar s(z)}$ is
the conjugation of $\rho$ by the homothety $\bar s(z)$. Then, by
Section~\ref{subsec:conjugation}, $\rho^{\bar s(z)}$ is
associated with data $(h(z)=z^{-2}, c, b(z)- \frac{\bar
s'(z)}{\bar s(z)})$. Hence, for $\bar s(z)$ satisfying $\frac{\bar
s'(z)}{\bar s(z)}=\sum\limits_{i\geq 0} b_{2i+1} z^{2i+1}$, we may assume
that $b_{2i+1}=0$. The same argument, for the conjugation by the
homothety $z^k$ with $k= b_{-1}-\frac32\in{\mathbb Z}$, shows that we may
assume $b_{-1}=\frac32$. That is, the desired series $s(z)$ equals
$z^k\bar s(z)$.
\end{Example}

\begin{Corollary}\label{cor:LieAlgLift}
Let $U\subset {\mathbb C}((z))$ be a subspace with $A_U={\mathbb C}[h(z)]$ and
stable under an action $({\mathbb C}((z)),\rho)$ defined by a triple
$(h(z),c,b(z))$.

There exists a semilinear transformation $g\in \operatorname{SGl}({\mathbb C}((z)))$ such
that $\widetilde A_{g(U)}^1$, the central extension of
$A_{g(U)}^1$, is split iff $b_{-1}\in\frac12{\mathbb Z}\setminus{\mathbb Z}$.
\end{Corollary}

\begin{proof}
The conclusion
follows from Theorem~\ref{thm:HomModAutxC((z))} and
Proposition~\ref{prop:LieAlgLift}.
\end{proof}

\begin{Remark}
The condition on the residue was also obtained, from another point of view, by Schwarz, see \cite[equation~(29)]{Schwarz}.
\end{Remark}

Let us say a word on the relevance of this corollary. Observe that $A_U^1$ acts on ${\mathbb C}((z))$ while its central extension, $\widetilde A_U^1$,  acts naturally on the fermionic Fock space
$H^0(\operatorname{Gr}(V),\operatorname{Det}_V^*)\simeq
\Lambda^{\frac{\infty}{2}}V$. Thus, in order to apply the bosonization isomorphism to relate subspaces and $\tau$-functions (as it was explained in the introduction), it is needed that operators (of our algebra) on ${\mathbb C}((z))$ can be lifted to operators (of the central extension) in a compatible way; that is, that the central extension~$\widetilde A_U^1$ is split and, thus, the restriction of the bosonization morphism to it is a Lie algebra homomorphism. Further, having in mind that~$U$ is invariant under the action of~$\widetilde A_U^1$, it follows that each element of~$\widetilde A_U^1$ def\/ines an homothety of the stalk of the determinant bundle at~$U$,~$\operatorname{Det}^*_U$. Thus, in order to def\/ine a canonical section of $\pi\colon\widetilde A_U^1\to  A_U^1$ one maps an element $a\in A_U^1$ to the element of $\pi^{-1}(a)$ such that induces the identity at~$\operatorname{Det}^*_U$. We will come again to this point in Section~\ref{subsec:taufunctions}.

\subsection{Stable subspaces}

In this subsection we aim to construct explicitly a subspace
fulf\/illing our requirements; namely, invariance under the action
and under the homothety~$z^{-2}$. Because of this fact and of Theorem~\ref{thm:stab=C[h]orC}, we shall assume, henceforth, that
$h(z)=z^{-2}$.

A naive candidate would be the ${\mathbb C}[h(z)]$-module generated by $1$ under
the action of~$\rho(L_{-1})$. Nevertheless, we shall need to
consider a conjugate of it (Section~\ref{subsec:conjugation}). For this, we shall choose a solution of the Airy
equation and decompose~$b(z)$ in a suitable way. Let us be more
precise.

First, we choose $w(z)$, a formal solution of the Airy equation
    \begin{gather}\label{eq:Airy}
    w''(z)   +  \frac12 S(h(z)) w(z)  =  0,
    \end{gather}
where $S$ denotes the Schwarzian derivative; that is
\begin{gather*}
    S(h) :=
    \frac{h'''(z)}{h'(z)}-\frac32\left(\frac{h''(z)}{h'(z)}\right)^2.
\end{gather*}
It is a straightforward check that $w(z)$ satisf\/ies the Airy
equation if\/f $f(z):=\frac{w'(z)}{w(z)}$ satisf\/ies the Riccati
equation
         \begin{gather}\label{eq:Riccati}
         f(z)^2 + f'(z) + \frac12 S(h(z))  =  0.
         \end{gather}
Note, in particular, that $h'(z)^{-1/2}$ satisf\/ies the equation~(\ref{eq:Airy}); or, equivalently, $\operatorname{d}\log (h'(z)^{-1/2}) =
-\frac12\frac{h''(z)}{h'(z)}$ satisf\/ies
equation~(\ref{eq:Riccati}).

Recalling from \cite[Chapters~6 and~9]{Laine} the basic properties
of the solutions of the Airy and Riccati equations, we know that
in our situation the solutions of equation~(\ref{eq:Riccati}) are
meromorphic; i.e., $\frac{w'(z)}{w(z)}\in {\mathbb C}((z))$. Thus, it makes
sense to conjugate a given action by $w(z)$ (see
Section~\ref{subsec:conjugation}).

Furthermore, if $w(z)$ and $h'(z)^{-\frac12}$ are linearly
independent, then the fact that they are solutions of the Airy
equation implies that the Schwarzian derivative of $w(z)
h'(z)^{\frac12}$ coincide with that of $h(z)$, and they therefore
dif\/fer by a M\"{o}bius transformation; that is
    \begin{gather}\label{eq:wh=rational}
    w(z)  h'(z)^{\frac12}  =
    \frac{\alpha     h(z)+\beta}{\gamma h(z)+\delta}
    \qquad\text{for some} \ \
        \begin{pmatrix} \alpha & \beta
        \\ \gamma & \delta
        \end{pmatrix}
    \in \operatorname{PGL}(2,{\mathbb C}).
    \end{gather}

Second, given a f\/irst-order dif\/ferential operator $P$ we search
for a formal expression $v(z)$ such that
\begin{gather*}
    P^v -u(h(z))   =  -\frac{1}{h'(z)}\partial_z
    -\frac12\frac{h''(z)}{h'(z)},
\end{gather*}
where $u(x)\in{\mathbb C}[x]$ (the reason is that the square of this
operator has no term in $\partial_z$). More concretely, for the
case  $P=-\frac{1}{h'(z)}\partial_z +\frac{b(z)}{h'(z)}$, we
express $b(z)$ w.r.t.\ the decomposition
\begin{gather*}
    {\mathbb C}((z))   \simeq   {\mathbb C}[h(z)]h'(z) \oplus
    \big({\mathbb C}[h(z)]+{\mathbb C}[[z]]z^{-1} \big)
    \simeq
    {\mathbb C}\big[z^{-2}\big]z^{-3}\oplus \big( {\mathbb C}\big[z^{-2}\big] + {\mathbb C}[[z]]z^{-1}\big),
\end{gather*}
since $h(z)=z^{-2}$ and $h'(z)=-2 z^{-3}$. That is, let us
consider the unique elements $u(h(z))$, $v(z)$, where $u(x)$ is a
polynomial and  $v(z)$ is a formal expression with
$\frac{v'(z)}{v(z)}\in {\mathbb C}[h(z)]+z^{-1}{\mathbb C}[[z]]$, such that
    \begin{gather}\label{eq:b=sum}
    b(z)  =
    u(h(z))h'(z)
     +
    \left( \frac{v'(z)}{v(z)} - \frac12{h''(z)}\right).
    \end{gather}
Equivalently,  $v(z)$ is def\/ined by the formal expression
\begin{gather*}
    v(z) :=
    \exp\int \left( b(z)
        - u(h(z))h'(z)
        + \frac12 h''(z)\right)\operatorname{d}z  .
\end{gather*}

\begin{Lemma}\label{lem:rho^2NEW}
Given $P=-\frac{1}{h'(z)}\partial_z+\frac{b(z)}{h'(z)}$, let
$w(z)$, $u(h(z))$, $v(z)$ as above. It then holds that
\begin{gather*}
    \big(
    P^2   -  2u(h(z)) P
    +  \big(u'(h(z))+u(h(z))^2\big)\big)
    (1\otimes w(z) v(z))
    =  0.
\end{gather*}
\end{Lemma}

\begin{proof}
Note that the l.h.s.\ in the statement is rewritten as
    \begin{gather*}
    \big(    P-  u(h(z)   )\big)^2  (1\otimes w(z) v(z))
    =
     \left(\left(    P-u(h(z))-\frac{1}{h'(z)}\frac{v'(z)}{v(z)}\right)^2
    (1\otimes w(z))\right) v(z)\\
\hphantom{\big(    P-  u(h(z)   )\big)^2  (1\otimes w(z) v(z))}{}
    =  \left(\left(-\frac{1}{h'(z)}\partial_z -
    \frac12\frac{h''(z)}{h'(z)}\right)^2 (1\otimes w(z))\right)  v(z)\\
\hphantom{\big(    P-  u(h(z)   )\big)^2  (1\otimes w(z) v(z))}{}
    =  \frac{1}{h'(z)^2}\left(\left(\partial_z^2 +
    \frac12 S(h(z))\right) (1\otimes w(z))\right) v(z).
    \end{gather*}

In order to see that the last expression vanishes, note that
    \begin{gather*}
    \partial_z^2(1\otimes w(z)) =
    \partial_z\big( \partial_z(1\otimes w(z))\big)  =
    \partial_z\left(\left(\partial_z+\frac{w'(z)}{w(z)} \right)(1)\otimes
    w(z)\right)
  \\
  \hphantom{\partial_z^2(1\otimes w(z))}{}
    =
    \partial_z\left(\frac{w'(z)}{w(z)}\otimes w(z)\right)
 =
    \left( \left(\frac{w'(z)}{w(z)}\right)' +
    \left(\frac{w'(z)}{w(z)}\right)^2\right)\otimes w(z)
     \\
  \hphantom{\partial_z^2(1\otimes w(z))}{}
    =  -\frac12 S(h(z))\otimes w(z),
    \end{gather*}
where the last equality comes from the fact that
$\frac{w'(z)}{w(z)}$ solves the Riccati
equation~(\ref{eq:Riccati}).
\end{proof}

\begin{Remark}
In \cite{KacSchwarz} the authors are able to solve the
second-order dif\/ferential equation $\big(\frac32\bar z+
\frac1{2\bar z}\partial_{\bar z}-\frac1{4\bar z^2}\big)^2\phi(\bar
z) =\bar z^2\phi(\bar z)$ (their $\bar z$ variable and our $z$
variable are related by $\bar z=(\frac13)^{\frac13}z^{-1}$) by
the substitution $\phi(\bar z)= \bar z^{1/2}\exp(\frac23{\bar
z}^{-3})\psi(\bar z)$ where $\psi(\bar z)$ is a solution of the
Airy equation. However, this makes sense since they show that
$\phi(\bar z)$ has an asymptotic expansion in ${\mathbb C}[[\bar z^{-1}]]$.
Observe that the previous Lemma can be thought of as an abstract
formalization of this substitution.
\end{Remark}

\begin{Theorem}[existence]\label{thm:ExistenceU(v)}
Let $w(z)$ be a solution of the Airy equation~\eqref{eq:Airy} linearly independent with $h'(z)^{-\frac12}$. Let $({\mathbb C}((z)),\rho)$ be defined by $(h(z)=z^{-2},c,b(z))$ and~let $v(z)$ be a formal function such that equation~\eqref{eq:b=sum} is fulfilled. Let $V^{wv}$ be the ${\mathbb C}((z))$-vector space ${\mathbb C}((z))\otimes w(z) v(z) $ with the conjugated action~$\rho^{wv}$.

It then holds that the ${\mathbb C}$-vector subspace of $V^{wv}$
\begin{gather*}
        {\mathcal U}(w)
          :=\langle 1\otimes w(z) v(z) , \rho^{wv}(L_{-1})(1\otimes w(z)  v(z))\rangle \otimes_{{\mathbb C}} {\mathbb C}[h(z)]
\end{gather*}
is ${\mathcal W}^+$-stable, it is a ${\mathbb C}[h(z)]$-module of rank $2$  and it   belongs to $\operatorname{Gr}(V^{wv})$.
\end{Theorem}

\begin{proof}
Let us denote $P:=\rho^{wv}(L_{-1})$. Lemma~\ref{lem:rho^2NEW} implies that ${\mathcal U}(w)$ is a $P$-stable ${\mathbb C}[h(z)]$-module. The ${\mathcal
W}^+$-stability follows from those facts and from the following
relations
     \begin{gather*}
     \rho^{wv}(L_i)  =   h(z)^i(h(z)\rho^{wv}( L_{-1})-(i+1)c),
     \\
     \rho^{wv}(L_{-1})(p(h(z))\otimes v(z))  =    p(h(z)) \rho^{wv}(L_{-1})(1\otimes v(z)) - p'(h(z))\otimes v(z).
     \end{gather*}

Let us prove that $P(1\otimes w(z) v(z)) \notin {\mathcal U}(w)\otimes_{{\mathbb C}[h(z)]} {\mathbb C}((h(z)^{-1}))$. Observe that
\begin{gather*}
    \frac{1}{h'(z)}\frac{w'(z)}{w(z)}  =
    \frac{1}{h'(z)}\operatorname{d}\log \big(w(z) h'(z)^{\frac12}\big) - \frac12\frac{h''(z)}{h'(z)^2}
    \in {\mathbb C}[[z^2]]
\end{gather*}
by equation~(\ref{eq:wh=rational}) (recall that $h(z)=z^{-2}$).  Computing how $P$ acts, we have
    \begin{gather*}
    P(f(z)  \otimes  w(z) v(z))  =
    \left(-\frac1{h'(z)}\partial_z +u(h(z)) - \frac12\frac{h''(z)}{h'(z)}-\frac{1}{h'(z)}\frac{w'(z)}{w(z)}\right)(f(z))\otimes w(z) v(z)
    \end{gather*}
and note that the term $\frac12\frac{h''(z)}{h'(z)}$ on the r.h.s.\ shifts the order by an odd number while all the other terms shift it by an even number. Hence, ${\mathcal U}(w)$ is a free ${\mathbb C}[h(z)]$-module of rank~$2$.

Finally, in order to prove that~${\mathcal U}(w)$  lies in the Sato Grassmannian, where we are considering $V^{wv}_+:= {\mathbb C}[[z]]\otimes w(z) v(z)$, one has to show the following two conditions; namely,
    \begin{gather*}
    \dim_{{\mathbb C}} \big({\mathbb C}[[z]]\otimes w(z) v(z)\cap  {\mathcal U}(w)\big)   < \infty,
\\
    \dim_{{\mathbb C}} {\mathbb C}((z))\otimes w(z) v(z)/
    \big({\mathbb C}[[z]]\otimes w(z) v(z) + {\mathcal U}(w)\big)   < \infty.
    \end{gather*}
Bearing in mind that $u(h(z)) - \frac12\frac{h''(z)}{h'(z)}-\frac{1}{h'(z)}\frac{w'(z)}{w(z)}$ does not belong to ${\mathbb C}((z^2))$, both conditions follow easily from the previous claims.
\end{proof}

The above constructed subspace depends clearly on the choice of a solution of the Airy equation. The following result studies what this dependence looks like.

\begin{Proposition}\label{prop:w1w2IsomGr}
Let $({\mathbb C}((z)),\rho)$ be an action of ${\mathcal W}^+$ defined by
the data $\{h(z)=z^{-2},c,b(z)\}$. Let  $w_1(z)$, $w_2(z)$ be two solutions
of~\eqref{eq:Airy}.

Then, up to ${\mathbb C}^*$,  there is a unique isomorphism of  ${\mathbb C}((z))$-vector spaces
$V^{w_1}\overset{\sim}\to V^{w_2}$ which is compatible w.r.t.\ the actions of the conjugated actions~$\rho^{w_1}$ and~$\rho^{w_2}$.
\end{Proposition}

\begin{proof}
We begin by constructing one isomorphism; we shall then prove the
uniqueness.

Let us consider $(V^{w_i},\rho^{w_i})$ as the conjugated action by
$w_i(z)$ (for $i=1,2$); that is, the ${\mathbb C}((z))$-vector space
$V^{w_i}$ is given by ${\mathbb C}((z))\otimes_{{\mathbb C}}{\mathbb C} w_i(z)$ and
$\rho^{w_i}$ by equation~(\ref{eq:conjugationRhoVequal}).

From
\cite[Chapter 6]{Laine}, we know that the fact that $w_1$, $w_2$
solve~(\ref{eq:Airy}) yields
\begin{gather*}
    S\left(\frac{w_1(z)}{w_2(z)}\right) =  S(h)
\end{gather*}
and that, therefore, there exists ${
\begin{pmatrix} \alpha & \beta \\ \gamma & \delta \end{pmatrix}
\in \operatorname{PGL}(2,{\mathbb C}) }$ such that
    \begin{gather}\label{eq:v1v2PGL}
    \frac{w_1(z)}{w_2(z)}  =
    \frac{\alpha h(z)+\beta}{\gamma h(z)+\delta}.
    \end{gather}

Let us now check that the ${\mathbb C}((z))$-linear map
    \begin{gather}
    {\mathbb C}((z))\otimes_{{\mathbb C}}{\mathbb C} w_1(z)   \longrightarrow  {\mathbb C}((z))\otimes_{{\mathbb C}}{\mathbb C} w_2(z),
    \nonumber\\
    1\otimes w_1(z)  \mapsto
    \frac{\alpha h(z)+\beta}{\gamma h(z)+\delta} \otimes w_2(z)\label{eq:LemIsomConju}
    \end{gather}
gives rise to an isomorphism that is compatible with the actions
of $\rho^{w_1}$ on the l.h.s.\ and of $\rho^{w_2}$ on the r.h.s.;
that is, one has to show that
\begin{gather*}
    \left(\frac{\alpha h(z)+\beta}{\gamma h(z)+\delta}\right)
    \rho^{w_1}(L_k)(f(z)\otimes w_1(z))
    =
    \rho^{w_2}(L_k)\left(\left(\frac{\alpha h(z)+\beta}{\gamma h(z)+\delta}\right)
    f(z)\otimes w_2(z)\right).
\end{gather*}
We shall only prove the case $k=-1$, $f(z)=1$, since the general
case goes along the same lines.

First, taking logarithms and derivatives in
equation~(\ref{eq:v1v2PGL}), we obtain
    \begin{gather}\label{eq:identity}
    \frac{w'_1(z)}{w_1(z)}  =  \frac{w'_2(z)}{w_2(z)} +
    \left(\frac{\alpha h(z)+\beta}{\gamma h(z)+\delta}\right)^{-1}
    \partial_z \left(\frac{\alpha h(z)+\beta}{\gamma h(z)+\delta}\right).
    \end{gather}

On the one hand, one computes the image of
\begin{gather*}
    \rho^{w_1}(L_{-1})(1\otimes w_1(z))  =  \left(-\frac{1}{h'(z)}\frac{w'_1(z)}{w_1(z)} + \frac{b(z)}{h'(z)}\right)\otimes w_1(z)
\end{gather*}
by the map~(\ref{eq:LemIsomConju}) and one obtains
    \begin{gather}
 \frac{\alpha h(z)+\beta}{\gamma h(z)+\delta}
    \left(-\frac{1}{h'(z)} \frac{w'_1(z)}{w_1(z)} + \frac{b(z)}{h'(z)}\right)\otimes w_2(z)
 \nonumber\\
   \qquad =
    -\frac{1}{h'(z)}
    \left(\frac{\alpha h(z)+\beta}{\gamma h(z)+\delta} \frac{w'_2(z)}{w_2(z)}
    +
    \partial_z\left(\frac{\alpha h(z)+\beta}{\gamma h(z)+\delta}\right)-
    \frac{\alpha h(z)+\beta}{\gamma h(z)+\delta}  b(z)
    \right)\otimes w_2(z),\label{eq:imagerho1}
    \end{gather}
where we have used the identity~(\ref{eq:identity}).

On the other hand, one has
    \begin{gather*}
     \rho^{w_2}(L_{-1}) \left(\frac{\alpha h(z)+\beta}{\gamma h(z)+\delta}\otimes w_2(z)\right)\\
     \qquad{}
 =  \left(\frac{\alpha h(z)+\beta}{\gamma h(z)+\delta}\right)
    \left( -\frac{1}{h'(z)}  \frac{w'_2(z)}{w_2(z)} + \frac{b(z)}{h'(z)}\right)\otimes w_2(z)
      -   \frac{1}{h'(z)} \partial_z \left(\frac{\alpha h(z)+\beta}{\gamma h(z)+\delta}
    \right)\otimes w_2(z)
    \end{gather*}
and, since this expression coincides with
equation~(\ref{eq:imagerho1}), it follows that~(\ref{eq:LemIsomConju}) is an isomorphism compatible with the actions.

Let us denote by $\phi$ the isomorphism (\ref{eq:LemIsomConju})
and let $\psi\colon V^{w_1}\to V^{w_2}$ be another isomorphism
compatible with the actions. The statement will be proved if we
can show that~$\phi\circ\psi^{-1}$ belongs to~${\mathbb C}^*$.

Let $f(z)$ be def\/ined by $\psi(1\otimes w_1(z))=f(z)\otimes
w_2(z)$. Thus
\begin{gather*}
    \big(\phi\circ\psi^{-1}\big) (1\otimes w_2(z))
 =  f(z)^{-1}\frac{\alpha h(z)+\beta}{\gamma h(z)+\delta}  \otimes w_2(z)
\end{gather*}
is a ${\mathbb C}((z))$-linear automorphism of $V^{w_2}$ that is compatible
with the action of $\rho^{w_2}$; that is,
$(\phi\circ\psi^{-1})\circ \rho^{w_2} = \rho^{w_2}\circ
(\phi\circ\psi^{-1})$ and, bearing in mind
Remark~\ref{rem:LieBraLkh^j}, it follows that
\begin{gather*}
    \frac1{h'(z)}\partial_z\left(f(z)^{-1}\frac{\alpha h(z)+\beta}{\gamma h(z)+\delta}\right)
   =  0
\end{gather*}
and, hence, $f(z)=\lambda   \frac{\alpha h(z)+\beta}{\gamma
h(z)+\delta}$ for $\lambda\in{\mathbb C}^*$ and the statement follows.
\end{proof}

\begin{Remark}
Let $(V,\rho)$ be an action of ${\mathcal W}^+$ under the
conditions of Theorem~\ref{thm:ExistenceU(v)}. Let us choose
$\gamma\in {\mathbb C}$. Bearing in mind Theorem~\ref{thm:VirasoroSubAlg},
we consider the action def\/ined by
\begin{gather*}
    \rho_{\gamma}(L_{k}) :=
    \rho(L_{k})+(k+1)h(z)^k\gamma\qquad\forall\, k\geq -1
\end{gather*}
and note that ${\mathcal U}(w)$ is stable under the action of
${\mathcal W}^+$ via $\rho_{\gamma}$ for all $\gamma$; that is,
there exists a~$1$-parameter family of actions preserving the same
subspace~${\mathcal U}(w)$.
\end{Remark}

\begin{Proposition}\label{prop:SameUimpliesSameb}
Let  $(V,\rho_i)$ $(i=1,2)$ be actions of ${\mathcal W}^+$ defined
by the data  $h_i(z)$, $c_i$, $b_i(z)$ as in Theorem~{\rm \ref{thm:VirasoroSubAlg}}
and let us assume that ${\mathfrak v}(h_i(z))<0$.

If there exists a subspace $U\in \operatorname{Gr}(V)$ which is
$\rho_i({\mathcal W}^+)$-stable for $i=1,2$, then there exist
$\alpha\in{\mathbb C}^*$, $\beta\in {\mathbb C}$ such that:
    \begin{enumerate}[$(i)$]\itemsep=0pt
    \item $h_1(z)=\alpha h_2(z)+\beta$;
    \item $\rho_1(L_{-1}) -\alpha \rho_2(L_{-1}) =
    \frac{1}{h'_1(z)}( b_1(z)-b_2(z))  \in A_U$.
    \end{enumerate}

Conversely, if there exist $\alpha\in{\mathbb C}^*$, $\beta\in {\mathbb C}$ such that $(i)$ and $(ii)$ hold, then a subspace $U\in \operatorname{Gr}(V)$ is $\rho_1({\mathcal W}^+)$-stable if and only if it is $\rho_2({\mathcal W}^+)$-stable.
\end{Proposition}

\begin{proof}
First, note that Theorem~\ref{thm:stab=C[h]orC} implies that
\begin{gather*}
    {\mathbb C}[h_1(z)]  =  A_U  =  {\mathbb C}[h_2(z)],
\end{gather*}
such that there exist $\alpha\in{\mathbb C}^*$, $\beta\in {\mathbb C}$  satisfying
the f\/irst item.

Observe that $\sigma(\rho_1(L_{-1}))
=\sigma(\alpha\rho_2(L_{-1}))=-\frac1{h'_1(z)}$ and, thus
\begin{gather*}
    \rho_1(L_{-1}) -\alpha \rho_2(L_{-1})  \in  {\mathbb C}((z)).
\end{gather*}
Since $U$ is stable under the $ \rho_1(L_{-1}) -\alpha
\rho_2(L_{-1})$, the element $\rho_1(L_{-1}) -\alpha
\rho_2(L_{-1})$ must belong to the stabilizer,~$A_U$.

It follows from Theorem~\ref{thm:stab=C[h]orC} that
\begin{gather*}
    \rho_1(L_{-1}) -\alpha \rho_2(L_{-1})  \in   {\mathbb C}[h_1(z)]
\end{gather*}
and the f\/irst part of the statement follows.

The converse follows from
Theorems~\ref{thm:stab=C[h]orC} and~\ref{thm:VirasoroSubAlg}.
\end{proof}

\subsection[$\tau$-functions]{$\boldsymbol{\tau}$-functions}\label{subsec:taufunctions}

It is common to work with the second derivatives of the logarithm
of the $\tau$-function instead of with the $\tau$-function itself.
Indeed, if one looks at the KP hierarchy expressed as Hirota
bilinear equations, one realizes that multiplication by the
exponential of a linear function preserves the set of solutions of
the hierarchy. Accordingly the hierarchy can be equivalently
stated for the second derivatives of the logarithm of the
$\tau$-function. Also, recall~\cite[Lemma~3.8]{SegalWilson}, which
states that, for $U\in\operatorname{Gr}(V)$, the function  $\exp(\sum_i
a_it_i)\tau_U(t)$ is the $\tau$-function of the point
$\exp(-\sum_i a_i\frac{z^i}{i})U\in\operatorname{Gr}(V)$; that is,
$\partial_{t_i}\partial_{t_j}\log \tau_U(t)$ does not
vary on the orbit of $U\in\operatorname{Gr}(V)$ under the action of ${\mathbb C}[[z]]^*$
by homotheties. Thus, our following results will be focused on the
study of $\partial_{t_i}\partial_{t_j}\log \tau_U(t)$
instead of on $\tau_U(t)$. Finally, it is worth noticing that
these expressions also appear in other topics such as the study of
monodromy-preserving deformations, Frobenius manifolds, etc.

As it was mentioned in the introduction, a key point in our approach
relies on the equivalence of conditions on subspaces (e.g., $\tilde
L U\subseteq U$) and on their corresponding $\tau$-functions (i.e.,
\mbox{$\bar L\tau_U(t)=0$}) where the operators are related through the
bosonization isomorphism. Nevertheless, now we are dealing with a
set of operators that generate a Lie algebra; namely,~$A_U^1$,
rather than just one operator. On the other hand, as it was
explained in Section~\ref{subsec:CentralExt}, the bosonization
isomorphism involves the central extension~$\widetilde A_U^1$.
Thus, in order to make of the bosonization isomorphism a~morphism
of Lie algebras, we demand the central extension~$\widetilde
A_U^1$ to be split.

Summing up, from now on, we shall assume that $h(z)=z^{-2}$ and
$b_{-1}\in\frac12{\mathbb Z}\setminus {\mathbb Z}$ (see
Proposition~\ref{prop:LieAlgLift} and
Corollary~\ref{cor:LieAlgLift}).

For an action $({\mathbb C}((z)),\rho)$ let $\bar\rho$ denote the action
induced on ${\mathbb C}[[t_1,t_2,\ldots]]$ through the bosonization
isomorphism; that is, $\bar\rho(L_k):=B\circ \rho(L_k)\circ
B^{-1}$ as a dif\/ferential operator on ${\mathbb C}[[t_1,\ldots]]$. Recall
that if $\rho$ corresponds to a triple $(h(z)=z^{-2},c,b(z))$ (KdV case), one usually deals only with
functions of $t$ variables with odd subindices and, therefore, if
no confusion arises, the composition
    \begin{gather}\label{eq:BarRhoOddVariables}
    \xymatrix@C=15pt{
    {\mathbb C}[[t_1,t_3,\ldots]] \ar@{^(->}[r] &
     {\mathbb C}[[t_1,t_2,\ldots]] \ar@{->>}[r]^-{\bar\rho(L_k)} &
     {\mathbb C}[[t_1,t_2,\ldots]]/(t_2,t_4,\ldots) \ar[r]^-\sim &
     {\mathbb C}[[t_1,t_3,\ldots]]}
     \end{gather}
will be denoted by $\bar\rho(L_k)$ too. Hence, $\bar\rho$ can be
thought of as acting on ${\mathbb C}[[t_1,t_3,\ldots]]$.

\begin{Theorem}[independence of choices]\label{thm:indep}
Let $({\mathbb C}((z)),\rho)$ be an action of ${\mathcal W}^+$ such that $h(z)=z^{-2}$ and
$b_{-1}\in\frac12{\mathbb Z}\setminus {\mathbb Z}$. Let  $w_1(z)$, $w_2(z)$ be solutions of the Airy equation~\eqref{eq:Airy}.

If the hypothesis of Theorem~{\rm \ref{thm:ExistenceU(v)}} holds and
${\mathcal U}(w_i)\in \operatorname{Gr} V^{w_iv}$ is the point given  by that
theorem $(i=1,2)$, then
\begin{gather*}
    \partial_{t_{2i+1}}\partial_{t_{2j+1}}\log  \tau_{{\mathcal U}(w_1)}(t)
    =
    \partial_{t_{2i+1}}\partial_{t_{2j+1}}\log  \tau_{{\mathcal U}(w_2)}(t)
    \qquad
    \forall\, i,j\geq 0
\end{gather*}
and this function is a common
solution of the KdV hierarchy and of the Virasoro constraint
equations
\begin{gather*}
    \bar \rho^{w_iv}(L_k) \tau_{{\mathcal U}(w_i)}(t)  = 0 \qquad \forall\, k\geq -1   ,\quad i=1,2.
\end{gather*}
\end{Theorem}

\begin{proof}
Note that the isomorphism $V^{w_1v}\simeq V^{w_2v}$
provided by Proposition~\ref{prop:w1w2IsomGr} yields an isomorphism $\operatorname{Gr}(V^{w_1v})\simeq \operatorname{Gr}(V^{w_2v})$ which sends ${\mathcal U}(w_1)$
to $\frac{\alpha h(z)+\beta}{\gamma h(z)+\delta}{\mathcal U}(w_2)$, where
\begin{gather*}
    \frac{w_1(z)}{w_2(z)} =
    \frac{\alpha h(z)+\beta}{\gamma h(z)+\delta}
\end{gather*}
(see the proof of the proposition).

To begin with, let us assume that
${\mathfrak{v}}\big(\frac{w_1(z)}{w_2(z)}\big)=0$ and, thus,
$\frac{\alpha h(z)+\beta}{\gamma h(z)+\delta}$ can be  expanded as
a series in $h(z)^{-1}$, say $g(h(z)^{-1})\in
{\mathbb C}[[h(z)^{-1}]]^*={\mathbb C}[[z^2]]^*\subseteq {\mathbb C}[[z]]^* $.

Bearing in mind~\cite[Lemma~3.8]{SegalWilson}, we know that the
$\tau$-function of  $g(h(z)^{-1}){\mathcal U}(w_2)$ is equal to
the $\tau$-function of ${\mathcal U}(w_2)$ up to the exponential
of a linear function on the $t$ variables. Hence, the second
derivatives of their logarithms do coincide, i.e.,
\begin{gather*}
\partial_{t_{2i+1}}\partial_{t_{2j+1}}\log
\tau_{{\mathcal U}(w_2)}(t)  =
\partial_{t_{2i+1}}\partial_{t_{2j+1}}\log
\tau_{g{\mathcal U}(w_2)}(t)  .
\end{gather*}

On the other hand, since
${\mathfrak{v}}\big(\frac{w_1(z)}{w_2(z)}\big)=0$ the isomorphism of the Proposition~\ref{prop:w1w2IsomGr}
sends $V^{w_1v}_+:={\mathbb C}[[z]]\otimes_{{\mathbb C}}{\mathbb C} w_1(z)\otimes_{{\mathbb C}}{\mathbb C} v(z)$ to
$V^{w_2v}_+:={\mathbb C}[[z]]\otimes_{{\mathbb C}}{\mathbb C} w_2(z)\otimes_{{\mathbb C}}{\mathbb C} v(z)$ and, hence, it is
compatible with the construction of the $\tau$-functions because
$\tau$-functions are def\/ined as the determinant of the projection
onto $V^{w_1v}/V^{w_1v}_+$, resp.\ $V^{w_2v}/V^{w_2v}_+$
(see~\cite{SegalWilson} or~\cite[Def\/inition~5.6]{AMP}). Therefore,
the $\tau$-function of ${\mathcal U}(w_1)$  and that of its image,
$g(h(z)^{-1}){\mathcal U}(w_2)$, coincide. The claim is proved.

It remains to check the cases where
${\mathfrak{v}}\big(\frac{w_1(z)}{w_2(z)}\big)\neq 0$. Observe
that there are two possibilities: either it is $2$ or $-2$. Let us
assume that it is $2$; that is, its expansion lies in
$h(z)^{-1}{\mathbb C}[[h(z)^{-1}]]^*= z^2 {\mathbb C}[[z]]^*$ and the isomorphism
induced between the Grassmannians sends the connected component of
index $\lambda$ to the connected component of index $\lambda+2$.
Recalling that if the $\tau$-function of a point in
$\operatorname{Gr}^{\lambda}(V^{w_2v})$ is the determinant of the projection onto
$V^{w_2v}/z^{\lambda}V^{w_2v}_+$, the $\tau$-function of a point in
$\operatorname{Gr}^{\lambda+2}(V^{w_2v})$ is therefore the determinant of the
projection onto $V^{w_2v}/z^{\lambda +2} V^{w_2v}_+$. Now, we can proceed as
above.
\end{proof}

\begin{Theorem}[uniqueness]\label{thm:uniqueness}
With $h(z)=z^{-2}$. For $i=1,2$, let  $({\mathbb C}((z)),\rho_i)$  be an action of ${\mathcal W}^+$ and let $\tau_i(t)\in {\mathbb C}[[t_1,t_3,\ldots]]$ be a $\tau$-function for the KdV hierarchy verifying Virasoro constraints
\begin{gather*}
    \bar\rho_i(L_k)(\tau_i(t))=0.
\end{gather*}

The following three conditions are equivalent:
\begin{enumerate}[$(i)$]\itemsep=0pt
    \item $\partial_{t_{2i+1}}\partial_{t_{2j+1}}\log  \tau_1(t) =    \partial_{t_{2i+1}}\partial_{t_{2j+1}}\log  \tau_2(t)$ for all $i,j\geq 0$;
    \item $ \rho_1(L_{-1})- \rho_2(L_{-1})   \in {\mathbb C}((z^2))$;
    \item $\bar\rho_1=\bar\rho_2$ on ${\mathbb C}[[t_1,t_3,\ldots]]$.
    \end{enumerate}
\end{Theorem}

\begin{proof}
First, since $\tau_i(t)$ is a $\tau$-function for the KdV
hierarchy, the Sato theory implies that it is the $\tau$-function
of a subspace $U_i\in\operatorname{Gr}({\mathbb C}((z)))$, such that $z^{-2}U_i\subset
U_i$.

Assuming that item (i) holds, we have that $\tau_1(t)$ and
$\tau_2(t)$ dif\/fer by a factor that is the exponential of a linear
function in the $t$ variables (with odd subindices). Lemma~3.8
of~\cite{SegalWilson} implies that there are $\alpha_i\in{\mathbb C}$ such
that $\tau_1(t)$ coincides with the $\tau$-function of $f(z) U_2$,
where $f(z)$ is $\exp\big({-}\sum\limits_{i\geq 0}
\alpha_{2i+1}\frac{z^{2i+1}}{2i+1}\big)$.

Hence, $\tau_1(t)$ is annihilated by $\operatorname{Im}\bar\rho_1$ and by
$\operatorname{Im}\bar\rho_2^{f}$,  with $\rho_2^{f}:=f(z)\circ\rho_2\circ
f(z)^{-1}$. Bearing in mind
Proposition~\ref{prop:SameUimpliesSameb} and the assumption that
$h_1(z)=h_2(z)=h(z)=z^{-2}$, it follows that
\begin{gather*}
    \rho_2(L_{-1})-\rho_2^{f} (L_{-1}) \in {\mathbb C}\big[z^{-2}\big].
\end{gather*}
Recalling the
explicit expressions (\ref{eq:rhoLk-hbc}),  we have{\samepage
    \begin{gather*}
    \rho_1(L_{-1})-   \rho_2^{f}(L_{-1})
     =  \frac{1}{h_1'(z)}\left( b_1(z)- \left(b_2(z)-
    \frac{\partial_z \exp\Big(\sum\limits_{i\geq 0} \alpha_{2i+1}\frac{z^{2i+1}}{2i+1}\Big)}
    {\exp\Big(\sum\limits_{i\geq 0} \alpha_{2i+1}\frac{z^{2i+1}}{2i+1}\Big)}\right)\right)
     \\
   \hphantom{\rho_1(L_{-1})-   \rho_2^{f}(L_{-1})}{}   =
    \frac{1}{h'(z)}\left( b_1(z)-b_2(z)+\sum_{i\geq 0} \alpha_{2i+1}z^{2i}\right),
    \end{gather*}
which lies in $ {\mathbb C}[h(z)]= {\mathbb C}[z^{-2}]$, and item (ii) is proved.}

Recall that the action of  $z^{-2m}$ (resp.\ $z^{2m}$) for $m>0$ on
the fermionic Fock space corresponds to the action of
$\partial_{t_{2m}}$ (resp. $2m t_{2m}$) on  the bosonic Fock space,
${\mathbb C}[[t_1,t_2,\ldots]]$ (this correspondence will be discussed in Section~\ref{sec:CFT}). Hence the operator $\bar
\rho(L_{-1})-\bar\rho_2 (L_{-1})$ acts on ${\mathbb C}[[t_1,t_2,\ldots]]$
as a linear combination of $\partial_{t_{2m}}$ and $2m t_{2m}$.
Considering their actions on ${\mathbb C}[[t_1,t_3,\ldots]]$ as in
(\ref{eq:BarRhoOddVariables}), we conclude item (iii).

Reversing the arguments, the converse implications (iii)$\implies$(ii)$\implies$(i)   follow easily.
\end{proof}

\begin{Remark}
Note that (iii) implies that for all $k\geq -1$ the systems of dif\/ferential
equations $\bar\rho_i(L_k)f(z)=0$ coincide for $i=1,2$. On the
other hand, if one replaces  the hypothesis $h_1(z)=h_2(z)=z^{-2}$ by
$h_1(z)=\alpha h_2(z)+\beta=z^{-2}$ for some $\alpha\in{\mathbb C}^*$ and
$\beta\in{\mathbb C}$, then (ii) must be replaced by ${\mathbb C}((z^2))+\operatorname{Im}\rho_1 =
{\mathbb C}((z^2))+\operatorname{Im}\rho_2$; and (iii) by $\operatorname{Im}\bar\rho_1=\operatorname{Im}\bar\rho_2$ on
${\mathbb C}[[t_1,t_3,\ldots]]$. These conditions are intimately related to
the f\/irst-order stabilizer (see Section~\ref{subsec:stab}).
\end{Remark}

The moral of this section is that we have been able to def\/ine
arrow $B$ of diagram~(\ref{eq:BigDiagramIntro}) of the Introduction and study some of
its properties. Summarizing, the previous Theorems show that: f\/irst, there
is an injection
    \begin{gather}\label{eq:DiagramMoral}
    \xymatrix@C=28pt@R=20pt{
    {\left\{\begin{gathered}
        % \text{functions }
        \tau(t)\in {\mathbb C}[[t_1,t_3,\ldots]] \\
        \text{satisfying KdV and string} \\
        \text{equation } \bar\rho(L_{-1})\tau(t)=0
    \end{gathered} \right\} }/ \sim
    \ar@<1ex>@{^(->}[r]^-{\bar A}
    &
    {\left\{\begin{gathered}
    \rho\in\operatorname{Hom}_{\text{Lie-alg}}\big({\mathcal W}^+,{\mathcal D}^1\big)
    \\
    \text{with }\sigma(\rho(L_{-1}))=-\frac12z^3
    \end{gathered}
    \right\} }/ \sim
    \ar@<1ex>@{-->}[l]^-{\bar B}
    }
    \end{gather}
where $\tau_1(t)$ and $\tau_2(t)$ are identif\/ied when the second derivatives of their logarithms coincide, and~$\rho_1$ and~$\rho_2$ are equivalent when $\rho_1(L_{-1})-\rho_2(L_{-1})\in{\mathbb C}((z^2))$. And second, that the arrow~$\bar B$ can be def\/ined on the subset of those~$\rho$ such that $b_{-1}\in\frac12{\mathbb Z}\setminus{\mathbb Z}$.

%%%%%%%%%%%%%%%%%%%%%%%%%%%%%%%%%%%%%%%%%%%%%
\section{Applications}\label{sec:ExamplesVirasoro}

\subsection{2D quantum gravity}\label{sec:CFT}

\looseness=1
Among the variety of topics in which Virasoro constraints and KdV
hierarchy appear together, we have chosen the case of case of 2D
quantum gravity within the framework of conformal f\/ield theory. We
do not  aim to review  the literature on CFT exhaustively, but to
illustrate how our results may help in the understanding of some
mathematical issues of the symmetries of that theory (e.g., the
partition function). Thus, our main references for this section
will be~\cite{DVV,Douglas,KacSchwarz,KSU,Kon}.  We shall see that
the Lie algebras considered by those authors correspond, in the
fermionic formulation, to representations of~${\mathcal W}^+$ on
$V={\mathbb C}((z))$ which f\/it into the approach of\/fered in
Section~\ref{sec:representationWitt}. Further, the data $h(z)$, $c$, $b(z)$ of
Theorem~\ref{thm:VirasoroSubAlg} will be written down explicitly.

Recall that the bosonization isomorphism, a.k.a.\ boson-fermion
correspondence (e.g., \cite[Lecture~5]{KacRa}, \cite[Section~2]{KSU}),
establishes an ${\mathbb C}$-linear isomorphism~$B$ between
$\Lambda^{\frac{\infty}{2}}{\mathbb C}((z))$ and ${\mathbb C}[[t_1,t_2,\ldots]]$
that preserves charge and degree. It is well known that $B$ can be
explicitly expressed in dif\/ferent f\/lavors such as Schur
polynomials, Matrix integrals or Laplace transforms. Therefore,
conjugating by $B$, any operator on the bosonic Fock space,
${\mathbb C}[[t_1,t_2,\ldots]]$, can be understood as an operator on the
fermionic Fock space, $\Lambda^{\frac{\infty}{2}}{\mathbb C}((z))$. It is
worth noticing that if we are given an algebra of operators on the
bosonic Fock space that is isomorphic to the Virasoro algebra (or to~${\mathcal W}^+$), then  the algebra generated by their
conjugations w.r.t.~$B$ is also isomorphic to the Virasoro algebra
(or to~${\mathcal W}^+$). As an example, we refer readers to
\cite[Table~1]{Kaz} where an explicit description on how the same
operator acts on both spaces is given.

In~\cite[Section~2.2]{KSU} the authors explain how the bosonic
formulation associated with the central charge $\mathbf c=1$
and the charge $0$ sector corresponds to the action of the f\/ields
    \begin{gather}
    \bar L_{-n}  =  \sum_{m=1}^{\infty}
    (n+m)t_{n+m}\partial_{t_m}+\frac12\sum_{m=1}^{n-1}
    m(n-m)t_{n-m} t_m,
    \qquad
    \bar L_0  = \sum_{m=1}^{\infty} m t_m\partial_{t_{m}},
    \nonumber\\
    \bar L_n =  \sum_{m=1}^{\infty} m
    t_m\partial_{t_{n+m}}+\frac12\sum_{m=1}^{n-1}\partial_{t_m}\partial_{t_{n-m}},
    \qquad
    {\mathbf c}=1\label{eq:KSU-VirOperators}
    \end{gather}
(where $n\geq 1$) in the Fock space ${\mathbb C}[[t_1,t_2,\ldots]]$.

Let us translate this picture into the fermionic formulation. It is well known the explicit correspondence between certain f\/irst-order dif\/ferential
operators on ${\mathbb C}((z))$ and certain dif\/ferential operators on the Fock space (see, for instance, \cite[Table~1]{Kaz}; see \cite{Givental} for an approach based on quantization). Indeed, if $a_m$ acts on ${\mathbb C}((z))$ by the homothety~$z^n$, then it acts on the Fock space~by
\begin{gather*}
    a_n :=
        \begin{cases}
        n t_n & \text{ for } n >0, \\
        0 & \text{ for } n=0, \\
        \partial_{t_{-n}} &\text{ for } n<0.
        \end{cases}
\end{gather*}

For $n\geq 1$, it holds that the action of the operator
$z^n(z\partial_z+\frac{1+n}2)$ acting on ${\mathbb C}((z))$ corresponds to the following operator on the bosonic Fock space
    \begin{gather*}
    \frac12 \sum_{m=-\infty}^{\infty} \colon a_m a_{n-m}\colon
     =
    \frac12 \sum_{m=1}^{n-1} a_m a_{n-m} + \frac12
    \sum_{m=n+1}^{\infty} a_m a_{n-m} +
    \frac12 \sum_{m=-1}^{-\infty} a_{n-m} a_{m}
     \\
\hphantom{\frac12 \sum_{m=-\infty}^{\infty} \colon a_m a_{n-m}\colon}{} =
    \frac12 \sum_{m=1}^{n-1} a_m a_{n-m} + \frac12
    \sum_{m=1}^{\infty} a_{n+m} a_{-m} +
    \frac12 \sum_{m=1}^{\infty} a_{n+m} a_{-m}
   \\
\hphantom{\frac12 \sum_{m=-\infty}^{\infty} \colon a_m a_{n-m}\colon}{}
 =
    \frac12 \sum_{m=1}^{n-1} m(n-m) t_m t_{n-m} +
    \sum_{m=1}^{\infty} (m+n) t_{n+m} \partial_{t_{m}}
    =  \bar L_{-n},
    \end{gather*}
where $\colon~~\colon$ denotes the normal ordering.
Analogously, the
action of $z^{-n}(z\partial_z+\frac{1-n}2)$ on ${\mathbb C}((z))$
corresponds to the action of
\begin{gather*}
    \frac12 \sum_{m=-\infty}^{\infty} \colon a_m a_{-n-m} \colon
 =  \bar L_{n}.
\end{gather*}

Hence, the action of ${\mathcal W}^+$ in $V={\mathbb C}((z))$ corresponding to $\{\bar L_n\vert n\geq -1\}$ is explicitly given by
    \begin{gather}\label{eq:KSU-ferminOper}
    \rho(L_n) :=   z^{-n}\left(z\partial_z+\frac{1-n}2\right)\qquad \forall\, n\geq -1.
    \end{gather}
Finally, let us compute $h(z)$, $b(z)$ and $c$ from the equations
    \begin{gather*}
  -\frac{h(z)^{n+1}}{h'(z)} = z^{-n+1},
\qquad
    -(n+1) c  h(z)^n + \frac{b(z)}{h'(z)} h(z)^{n+1} =
    \frac{1-n}2 z^{-n}
    \end{gather*}
and, therefore, $h(z)=z^{-1}$,  $c=\frac12$ and $b(z)= - z^{-1}$.

Let us now discuss how the KdV hierarchy shows up. Recall that,
for $m=2$ and the above data, Example~\ref{ex:VirSubalg} provides
a method to obtain another action with $h(z)=z^{-2}$, which
corresponds to the KdV. We shall see that this is indeed the case
for 2D gravity.

Douglas \cite{Douglas} proposed that the non-perturbative
partition function of two-dimensional gravity is the square of a
$\tau$-function for the KdV hierarchy that also satisf\/ies the
string equation.  In~\cite{DVV} the authors discussed further
consequences for the case of topological gravity but, again, the
constraints arised from the KdV hierarchy and from the string
equation. Let us see how. In their paper we f\/ind a Lie algebra,
that is isomorphic to ${\mathcal W}^+$, and that is generated by
the following operators (see~\cite[equation~(3.5)]{DVV})
    \begin{gather}
    \bar L'_{-1}  = \sum_{i=1}^{\infty} \left(i+\frac12\right){q_i}\partial_{{q_{i-1}}} + \frac18\lambda^{-2}q_0^2,
    \qquad
    \bar L'_0  = \sum_{i=0}^{\infty} \left(i+\frac12\right){q_i}\partial_{{q_i}} + \frac1{16},
   \nonumber \\
    \bar L'_{n}  = \sum_{i=0}^{\infty} \left(i+\frac12\right){q_i}\partial_{{q_{i+n}}} + \frac12\lambda^2\sum_{i=1}^n \partial_{{q_{i-1}}}\partial_{{q_{n-i}}}
    \qquad n\geq 1\label{eq:operartorsDVV}
    \end{gather}
acting on ${\mathbb Q}[{q_0},{q_1},\ldots]$. Note that this algebra coincides with the one considered by Givental~\cite[Section~3]{Givental} up to rescalling. Here, for the sake of simplicity, the constant term $\frac1{16}$ has been incorporated to $\bar L'_0$ so that $\{\bar L'_n\vert n\geq -1\}$ has Lie bracket $[\bar L'_i,\bar L'_j]=(i-j)\bar L'_{i+j}$ for $i,j\geq -1$.

Let us explain the relationship between the operators of
equation~(\ref{eq:KSU-VirOperators}) and those of
equation~(\ref{eq:operartorsDVV}). Indeed,  if one writes down how
the operators of equation~(\ref{eq:KSU-VirOperators}) act on the
subspace ${\mathbb C}[[t_1,t_3,t_5,\ldots]] \subset
{\mathbb C}[[t_1,t_2,t_3,\ldots]]$ (since we are dealing with
$\tau$-functions for the KdV hierarchy that depend only on $t$
variables with odd subindices), one has
    \begin{gather*}
    \bar L_{-2n}   =  \sum_{m=0}^{\infty}\!
    (2(n+m)+1)t_{2(n+m)+1}\partial_{t_{2m+1}}
      +   \frac12\!\sum_{m=0}^{n-1}\!
    (2m+1)(2(n-m)+1)t_{2(n-m)+1} t_{2m+1},
    \\
    \bar L_0  = \sum_{m=0}^{\infty} (2m+1) t_{2m+1}\partial_{t_{2m+1}},
    \\
    \bar L_{2n}  =  \sum_{m=0}^{\infty} (2m+1)
    t_{2m+1}\partial_{t_{2(n+m)+1}}+\frac12\sum_{m=0}^{n-1}\partial_{t_{2m+1}}\partial_{t_{2(n-m)+1}},\qquad
     {\bf c}   =   1,
    \end{gather*}
where $n\geq 1$. Letting ${q_m}:=\sqrt{2} \lambda t_{2m+1}$, we
observe that $\frac12 \bar L_{2n}+ \delta_{n,0}\frac1{16}$ coincide exactly with the operator~$\bar L_n'$ of~(\ref{eq:operartorsDVV}). Furthermore, the action of ${\mathcal W}^+$ on ${\mathbb C}((z))$
corresponding to the approach of \cite{DVV}; that, to the Lie algebra $\{\bar L'_n\,\vert\, n\geq -1\}$, is
    \begin{gather}\label{eq:Lk-KacSc}
    \rho'( L_{n}) :=   \frac12 \rho( L_{2n})  =
    \frac12 z^{-2n}\left(z\partial_z+\frac{1-2n}2\right)\qquad \forall\, n\geq -1,
    \end{gather}
which is the action def\/ined by the data $h(z)=z^{-2}$, $c=\frac12$
and $b(z)=- \frac32 z^{-1}$. Note that, as stated above, these
data can be obtained by the method of Example~\ref{ex:VirSubalg}
from $h(z)=z^{-1}$, $c=\frac12$ and $b(z)= - z^{-1}$ (def\/ining the
action~(\ref{eq:KSU-ferminOper})) with $m=2$.

\begin{Remark}
Observe that the
operators $\{\rho'( L_{n}) \}$  preserve the decomposition
${\mathbb C}((z))={\mathbb C}((z^2))$ $\oplus z{\mathbb C}((z^2))$. Indeed, recalling the
def\/initions $J(\bar z)=\sum\limits_{n=0}^{\infty} {q_n} \bar z^{n+1/2}$
\cite[formula~(2.15)]{DVV} and $L_m= \bar z^{-1/2}( \bar
z\frac{\partial}{\partial \bar z})^m \bar z^{-1/2}$
\cite[Section~3]{Givental} we observe that $\bar z^{1/2}$ plays a
relevant role in them (set $\bar z=z^2$ for a comparison with our
setup). It is worth pointing out that both
papers deal with KdV $\tau$-functions.
\end{Remark}

\begin{Remark}
In \cite{Kaz} the author aims to describe the relationship between
several equations for Hurwitz numbers and Hodge integrals. There, the f\/irst-order
dif\/ferential operators on the Fock space  are $\frac12 z^{-2n}\big(z\partial_z +
\frac{1-2n}2\big)$ for $ n\geq -1$ which generate a Virasoro algebra
\cite[Section~6]{Kaz}. From a straightforward computation it follows
that this action of ${\mathcal W}^+$ is def\/ined by $h(z)=z^{-2}$,
$c=\frac12$ and $b(z)=- \frac32 z^{-1}$ (applying
Theorem~\ref{thm:VirasoroSubAlg}) and it therefore coincides with
that of~\cite{DVV}.
\end{Remark}

At the time when the previously cited papers were published, Kac
and Schwarz~\cite{KacSchwarz} wondered about an explicit
description of a point in the Sato Grasmannian whose
$\tau$-function is the $\tau$-function arising in 2D gravity (see~\cite{DVV}).

First, let us note that the coordinate they used, $\bar z$, and
the one we are using are related by $\bar z=(\frac13)^{\frac13}
z^{-1}$. We shall rewrite their statements in terms of $z$ but,
for the sake of simplicity, let us start backwards.

The point in the Grassmannian is characterized by the fact that it
is the \emph{unique} point $U$ invariant under multiplication by
$z^{-2}$ (which is equivalent to the KdV) and under the action of
$A:= \frac32 \bar z + \frac12 \frac1{\bar z}\partial_{\bar z} -
\frac14\bar z^{-2} = -3^{\frac{2}{3}}\big({-}\frac12
z^{-1}+\frac12z^3\partial_z +\frac14 z^2\big)$ or, equivalently,
that $z^{-2}u(z)$, $A u(z)\in U$ for all $u(z)\in U$. It is worth
pointing out that ${\mathfrak v}(A u(z))={\mathfrak v}(u(z))-1$
for all $u(z)\neq 0$ and that
\begin{gather*}
    \big[A,z^{-2}\big]  =  3^{2/3}.
\end{gather*}

Hence, for a subspace $U\in \operatorname{Gr}{\mathbb C}((z))$ stable under the action of
$A$ and $z^{-2}$, let $u_0(z)\in U$ be the monic element such that
${\mathfrak v}(u_0(z))$ attains the maximum of $\{ {\mathfrak
v}(u(z))\vert u(z)\in U\}$. It then follows that $u_0(z)$, $A
u_0(z)$, $A^2 u_0(z)$ and $z^{-2}u_0(z)$ must be linearly dependent
over~${\mathbb C}$. After some computations the authors show that the
linear dependence implies that~$u_0(z)$ is a series in~$z^3$ with
${\mathfrak v}(u_0(z))=0$ (see \cite[equation~(8)]{KacSchwarz} for
the explicit expression of $u_0(z)$). Hence, it holds that
 \begin{gather}\label{eq:UKacSchwarz}
    U = \langle u_0(z), A u_0(z),\ldots, z^{2n}u_0(z),z^{2n} Au_0(z),\ldots\rangle .
    \end{gather}
Therefore, as claimed and proved in~\cite{KacSchwarz}, there is a
unique $U$ stable by $A$ and $z^{-2}$ and  it is explicitly given
as above.

We now wish to rephrase these arguments in terms of Witt algebras.
We wonder if there exists an action $({\mathbb C}((z)),\rho)$ of ${\mathcal
W}^+$ leaving $U$ (the subspace of
equation~(\ref{eq:UKacSchwarz})) stable. The fact that
$z^{-2}U\subset U$ implies that $h(z)=z^{-2}$. Since $U$ is
$A$-stable, it follows that $A$ belongs to the f\/irst-order
stabilizer of $U$ and, by the exact sequence (\ref{eq:A_U^1split})
and since $\sigma(A)=-3^{\frac{2}{3}}\frac12 z^3$, there exists a
polynomial $p(x)$ such that
\begin{gather*}
    A =  -3^{\frac{2}{3}} \rho(L_{-1}) +p(h(z)).
\end{gather*}
From where we obtain the equation
\begin{gather*}
    -\frac12z^{-1}+\frac14z^2
 =
    -\frac{1}{2}z^3 b(z) + p(h(z)).
\end{gather*}
We will consider the solution $b(z)=z^{-4} -\frac12z^{-1}$,
$p(h(z))=0$. Thus, we may assume that the action is def\/ined by
$h(z)=z^{-2}$, $c\in {\mathbb C}$ and $b(z) =z^{-4} -\frac12z^{-1}$.

In~\cite{Kon}, Kontsevich proved Witten's conjecture, which claims
that a generating function for intersection indices coincides with
the partition function in the standard matrix model theory. In
particular, it has to obey Virasoro constraints and the KdV
hierarchy. The action of these f\/ields on the space of fermions
corresponds to the following operators \cite[equation~(7.33)]{Dijkgraaf} (see also~\cite[equation~(2.59)]{Witten})
    \begin{gather*}
    \bar L'_{-1}  = -\frac12\partial_{q_0}+ \sum_{i=1}^{\infty} \left(i+\frac12\right){q_i}\partial_{q_{i-1}} + \frac14 q_0^2,
    \qquad
    \bar L'_0  = -\frac12 \partial_{q_1}+ \sum_{i=0}^{\infty} \left(i+\frac12\right){q_i}\partial_{{q_i}} + \frac1{16},
   \nonumber \\
    \bar L'_{n}  = -\frac12\partial_{q_{n+1}} +
    \sum_{i=0}^{\infty} \left(i+\frac12\right){q_i}\partial_{{q_{i+n}}} + \frac14\sum_{i=1}^n \partial_{{q_{i-1}}}\partial_{q_{n-i}},
    \qquad n\geq 1.%\label{eq:operartorsDijkgraff-WK}
    \end{gather*}
which have Lie bracket $[\bar L'_i,\bar L'_j]=(i-j)L_{i+j}$ for $i,j\geq -1$ and correspond to the following operators acting on ${\mathbb C}((z))$
\begin{gather*}
    \tilde L_n :=
    \frac12 z^{-2n}\left(z\partial_z+\frac{1-2n}2\right)
     - \frac12 z^{-(2n+3)},
  \qquad n\geq -1.
\end{gather*}

The Lie algebra generated by these operators is isomorphic to the Virasoro algebra and the Witten--Kontsevich $\tau$-function is annihilated by the
operators $\{\tilde L_n\}$ and is \emph{uniquely} characterized by
this property up to a scalar factor (see \cite[Section~3]{Givental}).

In our setting
    \begin{gather*}
      -\frac{h(z)^{n+1}}{h'(z)} = \frac12 z^{-2n+1},
    \qquad
      -(n+1) c   h(z)^n + \frac{b(z)}{h'(z)} h(z)^{n+1} =  \frac12 z^{-(2n+3)} -\frac{2n-1}4 z^{-2n}
    \end{gather*}
and, therefore, $h(z)=z^{-2}$,  $c=\frac12$ and $b(z)= z^{-4}-\frac32 z^{-1}$.

Notice that in the case of Kac--Schwarz~\cite{KacSchwarz} we
obtained $b(z)= z^{-4} - \frac12 z^{-1}$ while  we  now have
$b(z)= z^{-4} - \frac32 z^{-1}$. Let us explain why this is not
contradictory. The dif\/ference between both expressions is
$z^{-1}$, such that if the f\/irst case corresponds to the action on
a vector space $V$, the second one corresponds to the conjugated
action by $z\in{\mathbb C}((z))$ (see Theorem~\ref{thm:HomModAutxC((z))}).
Bearing in mind how $\tau$-functions on dif\/ferent connected components are computed, it follows that the $\tau$-functions of both
cases do coincide.

\begin{Remark}
It is worth mentioning the paper~\cite{Alexandrov}. Its Section~2 is concerned with  three instances of functions that solve simultaneously the KP hierarchy as well as Virasoro constraints; namely, Witten--Kontsevich, Hurwitz  and Hodge $\tau$-functions. For each case, he computes the so-called Kac--Schwarz operators; that is, those dif\/ferential operators preserving the point of the Grassmannian def\/ined by each of the above $\tau$-function. In the f\/irst case \cite[Section~2.3]{Alexandrov}, he obtains two Kac--Schwarz operators $a_{\rm KW}=\frac1z-z^3\frac{\partial}{\partial z}-\frac{z^2}2$, $b_{\rm KW}=z^{-2}$ which satisfy $[a_{\rm KW},b_{\rm KW}]_-=2$. Now, it is clear that the subalgebra generated by them coincide with~\eqref{eq:A1U}. Further~$A_U^1$ is a subalgebra of the Heisenberg--Virasoro algebra (see equation~\eqref{eq:LieBraLkh^j}). It is quite relevant that in the other two cases, the operators $a_{\rm Hurwitz}$, $b_{\rm Hurwitz}$ and $a_{\rm Hodge}$, $b_{\rm Hogde}$ are conjugated with $a_{\rm KW}$, $b_{\rm KW}$.
\end{Remark}

%%%%%%%%%%%%%%%%%%%%%%%%%%%%%%%%%%%%%%%%%%%%%%%%%%%%%%%%%%%
\subsection{Spectral curves}\label{subsec:spectral}

It is natural to wonder  whether the subspace of Theorem~\ref{thm:ExistenceU(v)}  can be described in terms of algebraic geometry via the Krichever morphism.

We continue with an action of ${\mathcal
W}^+$, $(V={\mathbb C}((z)),\rho)$, under the following assumptions:
 $h(z)=z^{-2}$; $v(z)$ and $u(x)$  satisfy equation~(\ref{eq:b=sum}) with $u''(x)\neq 0$; and the
hypotheses of Theorem~\ref{thm:ExistenceU(v)} hold. Let $U:={\mathcal U}(w)$ be the point of $\operatorname{Gr}(V^{wv})$ given by
Theorem~\ref{thm:ExistenceU(v)} for a solution $W$ of the Airy equation and let $A_U = {\mathbb C}[z^{-2}]$ be its
stabilizer.

Recall that by applying Krichever theory to the pair
$(A_U,U)$, one gets a rank $2$ vector bundle on the projective
line. Indeed, the projective line ${\mathbb P}_1$ is obtained by attaching
to ${\mathbb A}_1=\operatorname{Spec}A_U$ the point def\/ined by the
valuation ${\mathfrak v}$.
Let us set $\bar z:=h(z)$ as a coordinate in ${\mathbb P}_1\setminus \{{\mathbf 0}\}={\mathbb A}_1=\operatorname{Spec}A_U$.
Moreover, the $A_U$-module $U$, being a point in $\operatorname{Gr}(V^v)$, gives rise to a rank $2$ vector bundle on ${\mathbb P}_1$, ${\mathcal E}$.
Since $U$ is generated by $1\otimes v(z)$ and
$\rho^v(L_{-1})(1\otimes v(z))$ as a module over $A_U=
H^0({\mathbb P}_1\setminus \{{\mathbf 0}\},{\mathcal O}_{{\mathbb P}_1})$, let us introduce $T_{\rho}$, an
endomorphism of $U$ as an $A_U$-module,  by
  \begin{gather}
     T_{\rho}(1\otimes v(z))    :=   \rho^v(L_{-1})(1\otimes v(z)),
\nonumber
    \\
    T_{\rho}\big(\rho^v(L_{-1})(1\otimes v(z))\big)    :=   \rho^v(L_{-1})^2(1\otimes v(z))
\label{eq:Trho}.
    \end{gather}

In this case, $T_{\rho}$  is homogeneous of degree $d$ w.r.t.\ the f\/iltration $U\cap  \bar z^i {\mathbb C}[[z]]$. Thus, it yields a map
\begin{gather*}
    {\mathcal E}   \longrightarrow  {\mathcal E}\otimes {\mathcal O}_{{\mathbb P}_1}(d)
\end{gather*}
of ${\mathcal O}_{{\mathbb P}_1}$-modules or, equivalently, a map{\samepage
\begin{gather*}
    \operatorname{Sym}{\mathcal O}_{{\mathbb P}_1}(-d)   \longrightarrow  \operatorname{End}({\mathcal E})
\end{gather*}
of sheaves of ${\mathcal O}_{{\mathbb P}_1}$-algebras.}

Recalling the construction of the spectral curve~\cite{BNR},
these data def\/ine a 2:1-covering $\pi_{\rho}\colon X_{\rho}\to {\mathbb P}_1$
and a rank $1$ torsion free sheaf ${\mathcal L}_{\rho}$ on $X_{\rho}$ such that
$(\pi_{\rho})_*{\mathcal L}_{\rho} = {\mathcal E}$.

Let us look closer at this spectral curve. Recalling the explicit computations of the  proof of Lemma~\ref{lem:rho^2NEW},  it holds that  $T_{\rho}(\rho^v(L_{-1})(1\otimes v(z))) =-(u(\bar z)^2+u'(\bar z))\otimes v(z) + 2u(\bar z) \rho^v(L_{-1})(1\otimes v(z))$, and the restriction $T_{\rho}\vert_{{\mathbb P}_1\setminus\{{\mathbf 0}\}}$ has the associated matrix
\begin{gather*}
    \begin{pmatrix}
    0 & -u(\bar z)^2-u'(\bar z) \\ 1 & 2 u(\bar z)
    \end{pmatrix}
\end{gather*}
w.r.t.\ the basis $\{1\otimes v(z), \rho^v(L_{-1})(1\otimes
v(z))\}$ and thus the  characteristic polynomial of $T_{\rho}$ is
\begin{gather*}
    \operatorname{char}_{T_{\rho}}(y) =  y^2 - 2 u(\bar z) y
    + \big(u(\bar z)^2+ u'(\bar z)\big).
\end{gather*}
Hence, $X_{\rho}$ is an integral curve as long as $u'(x)$ is not a square. Furthermore, observe that
    \begin{gather}\label{eq:spectralcurve}
    X_{\rho}-\pi_{\rho}^{-1}(\{{\mathbf 0}\})
    \simeq
    \operatorname{Spec}{\mathbb C}[X,Y]/\big(Y^2 - 2 u(X) Y
    + \big(u(X)^2+ u'(X)\big)\big),
    \end{gather}
where the ring of the r.h.s.\ can be thought as a subalgebra of
$\operatorname{End}_{A_U}(U)$ identifying $X$ with the homothety of ratio
$\bar z=h(z)$ and $Y$ with the endomorphism~$T_{\rho}$. Finally,
the branch locus of $\pi_{\rho}$ consists of the zero locus of
$u'(\bar z)$ (recall that $\bar z$ is a coordinate in
${\mathbb P}_1\setminus \{{\mathbf 0}\}=\operatorname{Spec} A_U$) and the point~${\mathbf 0}$. The
singular points of~$X_{\rho}$ are precisely the preimages of the
locus of multiple roots of $u'(\bar z)=0$ together with the point
$\pi_{\rho}^{-1}(\{{\mathbf 0}\})$ (if $d>1$). Hence, $X_{\rho}$ is a
hyperelliptic curve and its genus can be easily expressed in terms
of~$u(\bar z)$.

The Hitchin map provided by the theory of spectral
curves takes the form
\begin{gather*}
    \rho \rightsquigarrow   \big( 2u(\bar z), -u(\bar z)^2-u'(\bar z)\big),
\end{gather*}
whose f\/ibers are essentially given by Jacobian varieties. Let us see this fact by studying~$b(z)$. Indeed,  if we express~$b(z)$ as in equation~(\ref{eq:b=sum})
\begin{gather*}
    b(z)  =
    u(h(z))h'(z)
     +
    \left( \frac{v'(z)}{v(z)} - \frac12{h''(z)}\right)
\end{gather*}
w.r.t.\ the splitting ${\mathbb C}((z))= {\mathbb C}[h(z)]h'(z)\oplus\big(
{\mathbb C}[h(z)]+z^{-1}{\mathbb C}[[z]]\big)$ we observe that varying the
coef\/f\/icients of $u$ may vary the spectral curve but not the
$\tau$-function (see Theorem~\ref{thm:uniqueness}) while the
variation of the second term, $\frac{v'(z)}{v(z)} - \frac12{h''(z)}$, leaves
$X_{\rho}$ f\/ixed but may change the sheaf~${\mathcal L}_{\rho}$ (this will be discussed in Section~\ref{subsec:UniverFamily}).
In other words, the variation of~$u(\bar z)$ would correspond to
monodromy-preserving tau functions and, hence, would be related to
Painlev\'{e} equations while the variation of the second term are KP
f\/lows.

\begin{Remark}
Note that
most of the examples given in \cite[Section~10]{EynardOrantin} are
based on curves that have the form of
equation~(\ref{eq:spectralcurve}) and that, in certain cases \cite[Sections~2 and~3]{EynardMarino}, the desired $\tau$-function,
satisfying Virasoro constraints and KdV hierarchy, can be
expressed in terms of the theta function of a hyperelliptic curve
(namely,~$X_\rho$). Thus, we shall investigate in a future paper  the relation between the
spectral curves and $\tau$-functions of our approach and the
spectral curves and $\tau$-functions of Eynard--Orantin~\cite{EynardOrantin,EynardOrantin2} since we expect that
both constructions are deeply connected.
\end{Remark}

\begin{Remark}
Very recently Schwarz has published \cite{Quantum}. In his Section~2 he shows that for a point~$U$ in the Grassmannian  invariant under $z^{q}$ and $\hbar \frac{d}{dz^p}+b(z)$ there are two dif\/ferential opera\-tors~$P$,~$Q$ s.t.\ $[P,Q]=\hbar$ and that they determine a \emph{quantum curve}. Observe that the study of this situation is almost equivalent, thanks to Proposition~\ref{prop:SameUimpliesSameb}, to the study of the case $\hbar =1$ (see~\cite{Schwarz}). Further, he considers a $z^q$-basis of~$U$; that is, vectors $v_0,\ldots,v_{q-1}$ such that $\{ z^{mq} v_n\}$ generate~$U$. It is easy to check that, for $q=2$, his~$U$ coincide with our~$U$ above and his equation~(3) with our~\eqref{eq:Trho}.
\end{Remark}

\subsection{Opers on the punctured disk}\label{subsec:oper}

Regarding the def\/inition and properties of the opers, which were
introduced by Drinfeld and Sokolov and generalized by Beilinson
and Drinfeld, we refer  interested readers to
\cite{FrLoop, Frenkel-Ben-Zvi}. For our own case, it suf\/f\/ices to
recall that a $\operatorname{Gl}(n)$-oper on the punctured disk,
$D^\times$, is a rank~$n$  vector bundle
$\mathcal E$ on  $D^\times$ equipped with a f\/lag $0\subset{\mathcal E}_1\subset\cdots\subset{\mathcal E}_{n-1}\subset{\mathcal E}_{n} = {\mathcal E}$ of subbundles  and a f\/lat  connection
$\nabla\colon {\mathcal E}\to {\mathcal E}\otimes\Omega_{D^\times}$ such that $\nabla({\mathcal E}_i)\subseteq {\mathcal E}_{i+1}\otimes \Omega_{D^\times}$ and the induced maps ${\mathcal E}_{i}/{\mathcal E}_{i-1}\to ({\mathcal E}_{i+1}/{\mathcal E}_i)\otimes\Omega_{D^\times}$ are a isomorphisms of ${\mathcal O}_{D^\times}$-bundles for all $i$ (transversality).

Let us see how we can associate an oper on the punctured disc $D^\times$ with certain
actions and conversely (see~\cite{Plaza-Oper} for a detailed exposition).

Let $(V={\mathbb C}((z)),\rho)$ be an action of ${\mathcal
W}^+$ such that $h(z)=z^{-2}$.  Let $U:={\mathcal U}(w)$ be the point of $\operatorname{Gr}(V^{wv})$ given by
Theorem~\ref{thm:ExistenceU(v)} and let $A_U = {\mathbb C}[z^{-2}]$ be its
stabilizer.

First, let ${\mathcal E}$  be the vector bundle on $\operatorname{Spec}{\mathbb C}((h(z)^{-1}))$ def\/ined by ${\mathbb C}((z))$. Hence, the subbundles ${\mathcal E}_i$, associated with
\begin{gather*}
    {\mathbb C}((h(z)^{-1}))\otimes_{{\mathbb C}}\langle 1, \rho(L_{-1})(1),\ldots, \rho(L_{-1})^{i-1}(1)\rangle
 \subseteq  {\mathbb C}((z)),\qquad i=1,2
\end{gather*}
def\/ine a f\/lag of  vector bundles $0\subset{\mathcal E}_1\subset {\mathcal E}_{2} = {\mathcal E}$.

Let us see that  $\mathcal E$ carries a connection. Indeed, let $\operatorname{d}\colon {\mathbb C}((h(z)^{-1}))\to
\Omega_{{\mathbb C}((h(z)^{-1}))/{\mathbb C}}$ be the dif\/ferential and consider the
${\mathbb C}$-linear map $ \rho (L_{-1})\otimes \operatorname{d}h - 1\otimes \operatorname{d}$:
    \begin{gather*}%\label{eq:connectionOper}
    {\mathbb C}((z))  \otimes_{{\mathbb C}} {\mathbb C}\big(\big(h(z)^{-1}\big)\big)    \longrightarrow  {\mathbb C}((z))  \otimes_{{\mathbb C}}\Omega_{{\mathbb C}((h(z)^{-1}))/{\mathbb C}},
    \\
    f\otimes a   \longmapsto
     \rho (L_{-1}) f \otimes a \operatorname{d}h - f\otimes \operatorname{d}a.
    \end{gather*}
One checks that when composing with the canonical map
\begin{gather*}
{\mathbb C}((z))  \otimes_{{\mathbb C}}\Omega_{{\mathbb C}((h(z)^{-1}))/{\mathbb C}} \to
{\mathbb C}((z))
\otimes_{{\mathbb C}((h(z)^{-1}))}\Omega_{{\mathbb C}((h(z)^{-1}))/{\mathbb C}}  ,
\end{gather*}
the images
of $f\otimes a$ and of $a f\otimes 1$ do coincide and, therefore, we
obtain a connection $\nabla\colon {\mathcal E}\to {\mathcal E}\otimes\Omega_{\operatorname{Spec}{\mathbb C}((h(z)^{-1}))}$ that satisf\/ies the transversality condition.

Let us see a converse; that is, given a $\operatorname{Gl}(2)$-oper on $D^{\times}$ one can def\/ine an action of ${\mathcal W}^+$. Let $({\mathcal E},\nabla)$ the $\operatorname{Gl}(2)$-oper on $D^{\times}=\operatorname{Spec}({\mathbb C}((h(z)^{-1})))$ where $h(z):=z^{-2}$. Let $\langle \,,\,\rangle $ be the pairing of dif\/ferentials with derivations. It can be shown that
\begin{gather*}
	\nabla_{D}(f) :=\langle \nabla f,D\rangle  \qquad \text{for} \ \ f\in {\mathbb C}((z))  , \quad  D\in \operatorname{Der}({\mathbb C}((h(z)^{-1})))
\end{gather*}
is a dif\/ferential operator of ${\mathbb C}((z))$ as a ${\mathbb C}((h(z)^{-1}))$-module. For $D_i:= -\frac{h(z)^{i+1}}{h'(z)}\frac{\partial}{\partial z}$, let us consider the linear map
    \begin{gather*}
    {\mathcal W}^+
   \overset{\rho}\longrightarrow
    {\mathcal D}^1_{{\mathbb C}((z))/{\mathbb C}}({\mathbb C}((z)), {\mathbb C}((z))),
    \\
    L_i  \mapsto   \rho(L_i):= \nabla_{D_i}.
    \end{gather*}
The fact that $\rho$ is a morphism of Lie algebras is derived from the f\/latness of $\nabla$ as follows
\begin{gather*}
	[\rho(L_i), \rho(L_j)]  =
	[ \nabla_{D_i}  ,  \nabla_{D_j} ]
	 =  \nabla_{[D_i,D_j]}  =
	\nabla_{(i-j)D_{i+j}}  =  (i-j) \rho(L_{i+j}).
\end{gather*}

However, the constructed oper  is def\/ined up to certain gauge transformations; namely, conjugation by $\operatorname{Aut}_{{\mathbb C}\text{-alg}}{\mathbb C}((\bar z))$ (in order to identify ${\mathbb C}((\bar z))$ and ${\mathbb C}((h(z)^{-1}))$) and by homotheties by~${\mathbb C}((z))$ (in order to identify~$V$ and~${\mathbb C}((z))$). The ef\/fects of these transformations have been studied in Section~\ref{subsec:conjugation}.

\begin{Remark}
Observe that, in the case $A_U={\mathbb C}[h(z)]$, the connection can be introduced in an alternative
way. Indeed, the map $h(z)^{n}\partial_h\mapsto L_{n-1}$, for
$n\geq 0$, provides a section of the canonical map $
{\mathcal D}^1_{A_U/{\mathbb C}}(U)\to \operatorname{Der}_{{\mathbb C}}(A_U)$ which, by
\cite[Section~1.1]{BeilinsonSchehtman}),  is an integrable connection on
${\mathcal E}$ on $\operatorname{Spec}{\mathbb C}((h(z)^{-1}))$.
\end{Remark}

\begin{Remark}
The techniques of \cite[Chapter~5]{Frenkel-Ben-Zvi} can be applied to the above results in order to associate to a a general action $(V,\rho)$ a $\operatorname{Gl}(n)$-oper on the abstract punctured disk $D^{\times}=\operatorname{Spec}{\mathbb C}((\bar z))$ (or, even, on an algebraic curve). On the other hand, recalling from \cite{KacPe,KacRa} the close relationship between vertex algebras and inf\/inite-dimensional representations of the Virasoro algebra, we expect to interpret the action of ${\mathcal W}^+$ on ${\mathbb C}[[t_1,t_2,\ldots]]$ in terms of vertex operators. This will help to understand our approach within Frenkel's framework of the geometric Langlands program~\cite{FrLoop}.
\end{Remark}

\subsection{Universal family}\label{subsec:UniverFamily}

Recent results on intersection theory are based on the so-called
\emph{topological recursion}
\cite{Kaz,LiuXu,MirzakhaniInvent,MirzakhaniJAMS,MulaseSafnuk,MulaseZhang}
which are formulae involving \emph{families} of $\tau$-functions
depending on an inf\/inite number of parameters such that the whole
family lies entirely on the space of functions satisfying KdV and
Virasoro constraints. One of these families already appeared in
Kontsevich's work~\cite[Section~3.4]{Kon}. It is worth mentioning the
existence of relevant $1$-parameter families; for instance, the  one \emph{connecting} Witten--Kontsevich and Mirzakhani theories~\cite{MulaseSafnuk}, and another one the Witten--Kontsevich partition function with the generating function of linear Hodge integrals def\/ined on the moduli space of stable curves~\cite{Kaz}.

It would be very interesting to generalize our methods to include the whole KP hierarchy and to interpret the family   \emph{connecting} the Witten--Kontsevich partition function with the Hurwitz partition function~\cite{MiMo} (see also~\cite{BoMar}). A recent result on this direction is in~\cite{Alexandrov} to be found.

In this subsection,  a natural procedure to
obtain the above mentioned families will be provided. Indeed, the proposed family will consist of a more general expression for $b(z)$ under the constraint that $b_{-1}\in\frac12{\mathbb Z}\setminus {\mathbb Z}$.

Let us consider a family of independent variables ${\mathbf s}:=(s_1,s_2,\ldots)$. For a sequence of non-negative integers, ${\mathbf m}:=(m_1,m_2,\ldots)$, with $m_i=0$ for all $i\gg 0$ def\/ine
\begin{gather*}
    \vert{\mathbf m}\vert :=  \sum_{i\geq 1} i m_i,
    \qquad
    \|{\mathbf m}\|  :=  \sum_{i\geq 1}  m_i,
    \qquad
    {\mathbf m}!  :=  \prod_{i\geq 1}  m_i!,
    \qquad
    {\mathbf s}^{\mathbf m} :=   \prod_{i\geq 1} s_i^{m_i}.
\end{gather*}

Based on Mulase--Safnuk's approach~\cite{MulaseSafnuk}, Liu--Xu considered the operators~\cite[equation~(9)]{LiuXu}
    \begin{gather*}
    \bar L'_n({\mathbf s}) :=
    -\frac12\sum_{{\mathbf m}}
    \frac{(-1)^{\|{\mathbf m}\|}}
        {{\mathbf m}! (2\vert{\mathbf m}\vert+1)!!}
    {\mathbf s}^{\mathbf m} \partial_{q_{\vert{\mathbf m}\vert +n+1}}
    +
    \sum_{i=0}^{\infty}\left(i+\frac12\right)  q_i \partial_{q_{i+n}}
    \\
\hphantom{\bar L'_n({\mathbf s}) :=}{}
 +
    \frac12\sum_{i=1}^n \partial_{q_{i-1}}\partial_{q_{n-i}}   +  \frac{{q_0}^2}{4}\delta_{n,-1}   +  \frac1{16}\delta_{n,0}
    \end{gather*}
for $n\geq -1$ (their exact expression corresponds to a rescaling by a double factorial). They showed that
\begin{gather*}
    [\bar L'_i({\mathbf s}) ,\bar L'_j({\mathbf s})]
    \,=\, (i-j) \bar L'_{i+j}({\mathbf s})
    \qquad\text{for} \ \ i,j\geq -1
\end{gather*}
and, therefore, they generate a family of Witt algebras depending on the parameters ${\mathbf s}$.

Observe that the operators $\bar L'_n(0)$, i.e., $\bar L'_n({\mathbf s})$ for ${\mathbf s}=0$, coincide with those of equation~(\ref{eq:operartorsDVV}) for $\lambda=2^{-\frac12}$ (up to rescaling of the variables $q_i$). That is, we have
\begin{gather*}
    \bar L'_n({\mathbf s})  =
    -\frac12\sum_{{\mathbf m}}
    \frac{(-1)^{\|{\mathbf m}\|}}
        {{\mathbf m}! (2\vert{\mathbf m}\vert+1)!!}
    {\mathbf s}^{\mathbf m} \partial_{q_{\vert{\mathbf m}\vert +n+1}}
 +   \bar L'_n(0).
\end{gather*}
Bearing in mind that the action induced by $\bar L'_n(0)$ on ${\mathbb C}((z))$ can be obtained by arguments analogous to those of the equation~(\ref{eq:Lk-KacSc}), we obtain the action corresponding to the above operators by replacing $q_m$ by $t_{2m+1}$ and $\partial_{t_{2m+1}}$ by $z^{-2m-1}$:
\begin{gather*}
    \rho'_{{\mathbf s}}(L_n) : = \, % &
    \frac12\sum_{{\mathbf m}}
    \frac{(-1)^{\|{\mathbf m}\|}}
        {{\mathbf m}! (2\vert{\mathbf m}\vert+1)!!}
    {\mathbf s}^{\mathbf m} z^{-2(\vert{\mathbf m}\vert +n)-3}
      +
    \frac12 z^{-2n}\left(z\partial_z+\frac{1-2n}2\right)\qquad \forall\, n\geq -1.
\end{gather*}
Since $\frac{h(z)^{n+1}}{h'(z)}=-\frac12z^{-2n+1}$, and regarding ${\mathbf s}$ as parameters, we obtain that the action $\rho'_{{\mathbf s}}$ is attached to the data $(h(z)=z^{-2}$, $c=\frac12$, $b_{{\mathbf s}}(z))$, where
\begin{gather*}
    b_{{\mathbf s}}(z) :=
    - \sum_{{\mathbf m}}
    \frac{(-1)^{\|{\mathbf m}\|}}
        {{\mathbf m}! (2\vert{\mathbf m}\vert+1)!!}
    {\mathbf s}^{\mathbf m} z^{-2\vert{\mathbf m}\vert-4}
       -  \frac32z^{-1}.
\end{gather*}

It must be pointed out that, bearing in mind Theorem~\ref{thm:HomModAutxC((z))}, Proposition~\ref{prop:LieAlgLift} and Theorem~\ref{thm:uniqueness}, we conclude that the action $\bar\rho'_{{\mathbf s}}$ (i.e., that induced by $\rho'_{{\mathbf s}}$ on ${\mathbb C}[[t_1,t_3,\ldots]]$) is the \emph{universal} action for the case of KdV (i.e., $h(z)=z^{-2}$). In particular, this agrees with the idea addressed in~\cite{MulaseSafnuk} that a certain $1$-parameter family, which would correspond to Eynard's \emph{spectral curve}, deforms the Witten--Kontsevich theory to other cases where the Virasoro also appears. Thus, our techniques provide an alternative proof of~\cite[Theorem~4.4]{MulaseSafnuk} (see also \cite[Theorem~2.1]{Kaz} and \cite[Theorem~4.1]{LiuXu}).

\begin{Theorem}
Let $\tau_{\mathbf s}(t)\in{\mathbb C}[[t_1,t_3,\ldots]]$ be the $\tau$-function associated to $\rho'_{\mathbf s}$. Then, $\tau_{\mathbf s}(t)$ sa\-tisfies the Virasoro constraints corresponding to operators $\bar L'_n({\mathbf s})$ above $($as in~Section~{\rm \ref{subsec:taufunctions})} and, moreover, it holds that
\begin{gather*}
    \tau_{\mathbf s}(t)
     =
    \tau_{0}(\tilde t),
\end{gather*}
where $\tilde t_{2+1}$ is equal to $t_i$ for $i=0,1$ and to $t_{2i+1}- \frac1{(2i+1)!!}\sum\limits_{\vert{\mathbf m}\vert=i-1 }\frac{(-1)^{\|{\mathbf m}\|}} {{\mathbf m}!} {\mathbf s}^{\mathbf m}$ for $i>1$.

Conversely, let $\tau(t)\in{\mathbb C}[[t_1,t_3,\ldots]]$ be a $\tau$-function for the KdV satisfying the Virasoro constraints. Then, there exist values of ${\mathbf s}$, say  ${\mathbf s}_0:=(s_1,s_2,\ldots)$ such that
\begin{gather*}
    \tau(t)
     =
    \tau_{{\mathbf s}_0}(t).
\end{gather*}
\end{Theorem}

\begin{proof}
It is enough to observe that under that change of variables, the operators $\bar L'_n({\mathbf s})$ in~$t_{2i+1}$ are transformed into the operators~$\bar L'_n(0)$ in~$\tilde t_{2i+1}$. Bearing in mind that, regarding~${\mathbf s}$ as para\-me\-ters, the action~$\rho'_{{\mathbf s}}$ is universal, the converse follows.
\end{proof}

Finally, note that the variation of ${\mathbf s}$ is equivalent to the variation of the line bundle constructed in Section~\ref{subsec:spectral}. Further, the previous approach can be also understood by replacing the base f\/ield ${\mathbb C}$ by the ring ${\mathbb C}[[s_1,s_2,\ldots]]$. Then, the resulting framework would make use of the relative Sato Grassmannian while the change of variables $t_{2+1}$ by $\tilde t_{2+1}$ corresponds to a trivialization of the relative Grassmannian~\cite{Plaza-RelGr}.

\subsection*{Acknowledgements}

I wish to thank Edward Frenkel, Esteban G\'{o}mez and, very specially, Motohico Mulase for stimulating discussions and to the referees for their remarks.

\pdfbookmark[1]{References}{ref}
\LastPageEnding

\end{document}